\documentclass[fontsize=11pt,twoside]{scrreprt}
\usepackage[square,numbers,nonamebreak]{natbib}
\usepackage{setspace}
\usepackage[english]{babel}
\usepackage{chngcntr}
\counterwithout{equation}{section}
\usepackage[hang]{footmisc}
 at 11truept
 at 11 truept

\usepackage[protrusion=true, expansion=true, final, kerning=true, spacing=true, factor=1100, stretch=10, shrink=10, babel, tracking=alltext, letterspace = 10]{microtype}
\usepackage[hmarginratio=1:1, bottom=2.5cm, top=2.5cm, inner=2.5cm, outer=2.5cm, paperwidth=17cm, paperheight=24cm]{geometry}

\usepackage[utf8]{inputenc}
\usepackage[sc,osf]{mathpazo}
\usepackage[euler-digits]{eulervm}
\usepackage[T1]{fontenc}
\usepackage{xcolor,graphicx,amsfonts,amsmath,ragged2e,titlesec,mathtools}

\usepackage[automark]{scrlayer-scrpage}
\cohead{\textnormal{Permutation groups, derangements and connections}}
\lehead{\pagemark}
\rohead{\pagemark}
\let\ceheadL\cehead
\renewcommand\cehead[1]{\ceheadL{\textnormal{#1}}}
\cehead{Scott Harper}
\lefoot{}
\rofoot{}

\newcommand\auths[1]{\large \textsc{#1}\setstretch{1.2}}
\newcommand\titl[1]{\center \linespread{1.1}\Large\textbf{ #1}\bigskip} 
\renewcommand\abstract[1]{
\begin{center}
{\textbf{Abstract}}
\end{center}
{
\linespread{1.1}\fontsize{9pt}{-10pt}\selectfont #1}}

\DeclareSymbolFont{operators}{\encodingdefault}{ppl}{m}{n}
\DeclareMathAlphabet{\mathbf}{\encodingdefault}{ppl}{bx}{n}
\DeclareMathAlphabet{\mathit}{\encodingdefault}{ppl}{m}{it}

\renewcommand{\thesection}{\arabic{section}}
 
\titleformat{\section}{\medskip\bigskip\normalfont\Large\bf}{\thesection}{0.5em}{}
\titleformat{\subsection}{\smallskip\bigskip\normalfont\large\bf}{\thesubsection}{0.5em}{}
\titleformat{\subsubsection}{\smallskip\bigskip\normalfont\large\bf}{}{0.5em}{}

\usepackage{amsthm}
\newtheoremstyle{dotless}{}{}{\itshape}{}{\bfseries}{}{1em}{}

\DeclareOldFontCommand{\rm}{\normalfont\rmfamily}{\mathrm}
\DeclareOldFontCommand{\sf}{\normalfont\sffamily}{\mathsf}
\DeclareOldFontCommand{\tt}{\normalfont\ttfamily}{\mathtt}
\DeclareOldFontCommand{\bf}{\normalfont\bfseries}{\mathbf}
\DeclareOldFontCommand{\it}{\normalfont\itshape}{\mathit}
\DeclareOldFontCommand{\sl}{\normalfont\slshape}{\@nomath\sl}
\DeclareOldFontCommand{\sc}{\normalfont\scshape}{\@nomath\sc}

\theoremstyle{dotless}
\newtheorem{theo}{Theorem}
\newtheorem{lem}[theo]{Lemma}
\newtheorem{cor}[theo]{Corollary}
\newtheorem{rem}[theo]{Remark}
\newtheorem{ex}[theo]{Example}
\newtheorem{que}[theo]{Question} 
\newtheorem{conj}[theo]{Conjecture} 
\newtheorem{exer}[theo]{Exercise} 

\counterwithout{theo}{section}
\numberwithin{theo}{section}

\usepackage{stmaryrd,multicol,tikz}

\newcommand{\Nat}{\mathbb{N}}
\newcommand{\Int}{\mathbb{Z}}
\newcommand{\Rat}{\mathbb{Q}}
\newcommand{\QQ}{\overline{\mathbb{Q}}}
\newcommand{\Real}{\mathbb{R}}
\newcommand{\Com}{\mathbb{C}}
\newcommand{\Ham}{\mathbb{H}}
\newcommand{\F}{\mathbb{F}}
\newcommand{\FF}{\overline{\F}}
\newcommand{\Primes}{\mathbb{P}}
\renewcommand{\P}{\mathcal{P}}

\renewcommand{\leq}{\le}
\renewcommand{\geq}{\ge}

\DeclareMathSymbol{\nmid}{\mathrel}{AMSb}{"2D}

\newcommand{\im}{\operatorname{im}}

\newcommand{\Aut}{\mathrm{Aut}}
\newcommand{\Gal}{\mathrm{Gal}}
\newcommand{\fix}{\mathrm{fix}}
\newcommand{\Sym}{\mathrm{Sym}}
\newcommand{\FSym}{\mathrm{FSym}}
\newcommand{\Alt}{\mathrm{Alt}}
\newcommand{\Frob}{\mathrm{Frob}}

\newcommand{\SL}{\mathrm{SL}}
\newcommand{\GL}{\mathrm{GL}}
\newcommand{\AGL}{\mathrm{AGL}}
\newcommand{\GaL}{\Gamma\mathrm{L}}
\newcommand{\AGaL}{\mathrm{A}\Gamma\mathrm{L}}
\newcommand{\PSL}{\mathrm{PSL}}
\newcommand{\PGL}{\mathrm{PGL}}
\newcommand{\Sp}{\mathrm{Sp}}
\newcommand{\PSp}{\mathrm{PSp}}
\newcommand{\PSU}{\mathrm{PSU}}
\newcommand{\POm}{\mathrm{P}\Omega}

\newcommand{\sectionquestion}[1]{%
  \begin{center}
  \begin{minipage}{0.8\textwidth}\centering\itshape #1\end{minipage}
  \end{center}\vspace{-2pt}}

\usepackage[hidelinks]{hyperref}

\begin{document}
\setheadsepline{1pt}[\color{black}]

\quad

\vspace{36.5mm}

\titl{Permutation groups, derangements \\ and connections}

\auths{Scott Harper}
\center{\phantom{blank}}
\thispagestyle{empty}\justify\noindent\setstretch{0.3}
\abstract{
What is the chance that a random shuffling of a deck of cards leaves no card in its original place? 
More precisely, what is the chance that a random permutation has no fixed points? 
A permutation with no fixed points is known as a derangement. 
Questions about derangements, like the ones above, go back centuries, but there are also several unsolved problems in this very active area of permutation group theory.

Based on a minicourse given at ``Groups and Algebras in Bicocca for Young Algebraists'' in June 2026, this survey provides a general introduction to permutation group theory, with a specific focus on derangements.
The survey features modern methods for studying permutation groups, recent developments and open problems concerning derangements and connections between this topic and other areas, such as combinatorics, graph theory, number theory, representation theory and topology.
}

\vspace{5.2mm}

\setstretch{1.1}
{\fontsize{11.5pt}{13.5pt}\selectfont


\section*{Introduction}

Permutation groups are where group theory began in the 19th century, and derangements (permutations with no fixed points) have been studied for even longer, going back at least to de~Montmort's 1708 book on games of chance.
That said, many open questions still remain about this deceptively simple notion, and progress on this topic shows how the general theory of permutation groups has developed over the years.
The simplicity of their definition also means that derangements arise in many areas of mathematics.
\clearpage

These themes are explored in this survey, which arose from a minicourse for early career researchers in algebra in June 2026.

\emph{Structure.}
The survey is divided into three sections.
Section 1 gives background on permutation groups (and can be omitted by anyone familiar with them), Section 2 discusses what is and is not known about derangements and Section 3 highlights twelve connections with other areas.

\emph{Coverage.}
The coverage is very selective, limited to finite transitive groups and focussing on two main themes: the proportion of derangements and derangements of prime-power order. 
These themes alone give rise to recent breakthroughs, open questions and connections with other areas.
However, much more could have been said: all discussion of the subgroup generated by derangements, width questions and the number of conjugacy classes of derangements have been omitted (see \cite[Section~1]{ref:BurnessFusari25} for more on these) as have the more combinatorial aspects of derangements (see \cite{ref:Cameron07} for more on these).
That said, attempts have been made to be more comprehensive with the connections, mentioning several that did not feature in the minicourse.

\emph{Reading order.}
To avoid all theory coming at the start and all the applications coming at the end, the minicourse proceeded roughly as follows, with the material that was omitted, or given without proof, in brackets:
\begin{itemize}
\item[(1)] 
derangements:
\ref{ss:derangements_history}, \ref{ss:permutations_basics}, \ref{ss:permutations_orbits}, \ref{ss:derangements_existence}, \ref{ss:connections_abstract}, \ref{ss:connections_toplogy},  (\ref{ss:connections_kronecker}), (\ref{ss:connections_intersecting})
\item[(2)]
proportion:
\ref{ss:permutations_tuples}, \ref{ss:derangements_proportion_1}, \ref{ss:connections_galois}, (\ref{ss:derangements_proportion_2}), (\ref{ss:connections_values})
\item[(3)]
prime-power order:
\ref{ss:permutations_equivalence}, \ref{ss:permutations_primitivity},  \ref{ss:derangements_prime_power}, \ref{ss:connections_brauer},  (\ref{ss:connections_character_vanishing}), (\ref{ss:connections_character_restriction}), \ref{ss:connections_isbell}
\item[(4)]
prime order: 
(\ref{ss:permutations_simple}), \ref{ss:derangements_prime}, \ref{ss:connections_polycirculant}.
\end{itemize} 

\emph{Further reading.}
Many references are given in the body of the text, but some general comments are in order.
The general theory in Section~\ref{s:permutations} is covered in more detail by Cameron's lovely book \cite[Sections~1.1--1.4 \& 1.8--1.12]{ref:Cameron99}; see also Wielandt's classical text~\cite{ref:Wielandt64} and Dixon and Mortimer's comprehensive graduate textbook~\cite{ref:DixonMortimer96}.
Sections~\ref{s:derangements} and~\ref{s:connections} benefited from a number of existing surveys including \cite{ref:ArezoomandAbdollahiSpiga19, ref:Cameron07, ref:Cameron11, ref:DiaconisFulmanGuralnick08, ref:Serre03} and the first chapters of \cite{ref:BubboloniSpigaWeigel24,ref:BurnessGiudici16}.

\clearpage

\section{Permutation groups}
\label{s:permutations}

The foundational work of Galois and Cauchy in group theory just concerned subgroups of symmetric groups, what we now call \emph{permutation groups} and what 19th-century writers called \emph{substitution groups}. 
Indeed, Lagrange's Theorem was first proved just for permutation groups, by Jordan in 1861 \cite{ref:Jordan61,ref:Roth01}. 
Cayley was one of the first to realise that groups could be studied through their key abstract properties, and
in 1854 he essentially proved what we now call Cayley's Theorem \cite{ref:Cayley54}, that every abstract group is isomorphic to a
permutation group.

\subsection{Basics}
\label{ss:permutations_basics}

Let $\Omega$ be a set.
The \emph{symmetric group}~$\Sym(\Omega)$ is the group of all permutations of~$\Omega$ under composition, and a \emph{permutation group} on~$\Omega$ is a subgroup of~$\Sym(\Omega)$.
When $\Omega = \{ 1, \dots, n \}$, we write $\Sym(\Omega) = \Sym(n)$ (and $\Alt(\Omega) = \Alt(n)$).
An \emph{action} of a group~$G$ on~$\Omega$ is a homomorphism $\rho \colon G \to \Sym(\Omega)$, and for $\omega \in \Omega$ and $g \in G$, we simply write $\omega g$ for $\omega (g \rho)$. 
For an action $\rho \colon G \to \Sym(\Omega)$, we write $G^\Omega$ for $\im(\rho)$, which is a permutation group, and we say that the action~$\rho$ is \emph{faithful} if $\ker(\rho)$ is trivial. 
The \emph{degree} of an action on~$\Omega$ is the cardinality of~$\Omega$.
Terms defined for actions also apply to the corresponding permutation group.
The set $\Omega$ will always be nonempty and usually finite.

Let a group $G$ act on a set~$\Omega$.
The \emph{orbit} of~$\omega \in \Omega$ is the subset 
\[ 
\omega G = \{ \omega g  \mid g \in G \},
\]
the \emph{stabiliser} of~$\omega$ is the subgroup 
\[
G_\omega = \{ g \in G \mid \omega g = \omega \}
\]
and the Orbit--Stabiliser Theorem tells us that $|G:G_\omega| = |\omega G|$, witnessed by the bijection
\[
G/G_\omega \to \omega G \qquad G_\omega x \mapsto \omega x.
\]
The action is \emph{transitive} if $\omega G = \Omega$ for all $\omega \in \Omega$ (i.e. there is a unique orbit), \emph{semiregular} if $G_\omega = 1$ for all $\omega \in \Omega$ (i.e. no nonidentity element fixes a point) and \emph{regular} if it is transitive and semiregular (so, by the Orbit--Stabiliser Theorem, $|G| = |\Omega|$).

\subsection{Orbit counting}
\label{ss:permutations_orbits}

We now give a tool for counting up to symmetry, which we will just use in a somewhat formal way.
For more on this technique and its generalisations, see \cite[Section~2.3]{ref:Cameron99} and \cite[Chapter~7]{ref:Cameron17}.

The key lemma is often referred to as Burnside's Lemma.
As Neumann \cite{ref:Neumann79} explains, this name, which seems to have arisen in the 1960s, is misleading.
Although Burnside states this lemma in his 1897 book \cite{ref:Burnside97}, and its 1911 second edition \cite{ref:Burnside11}, it was already known to Frobenius in 1887 \cite{ref:Frobenius87} and, to some extent, to Cauchy in 1845 \cite{ref:Cauchy45}.
Indeed, the first edition credits Frobenius \cite[119]{ref:Burnside97}.

In this lemma, and throughout, for $g \in \Sym(\Omega)$, we write
\[
\fix(g) = |\{ \omega \mid \omega g = \omega \}|.
\]

\begin{lem}[Orbit-Counting Lemma]
\label{lem:orbit-counting}
Let $\Omega$ be a finite set and let $G \leq \Sym(\Omega)$.
Then
\[
\text{the number of orbits of~$G$} = \frac{1}{|G|} \sum_{g \in G} \fix(g).
\]
\end{lem}

\begin{proof}
We double count the following set 
\[
X = \{ (\omega, g) \in \Omega \times G \mid \omega g = \omega \}
\]
as $|X| = \sum_{\omega \in \Omega} |G_\omega|$ and $|X| = \sum_{g \in G} \fix(g)$.
Equating these and using the Orbit--Stabiliser Theorem gives
\[
\frac{1}{|G|} \sum_{g \in G} \fix(g) 
= \frac{1}{|G|} \sum_{\omega \in \Omega} |G_\omega|
= \sum_{\omega \in \Omega} \frac{1}{|\omega G|},
\]
which is easily seen to be the number of orbits of~$G$.
\end{proof}

\subsection{Actions on tuples}
\label{ss:permutations_tuples}

Let $G$ be a group acting transitively on a set $\Omega$.
For $k \geq 1$, the group $G$ acts componentwise on the set $\Omega^k$ of $k$-tuples of elements of~$\Omega$ and the set $\Omega^{(k)}$ of $k$-tuples of distinct elements of~$\Omega$.

An \emph{orbital} is an orbit on $\Omega^2$, and the \emph{rank} of the action is the number of orbitals.
The \emph{diagonal orbital} is $\{ (\omega, \omega) \mid \omega \in \Omega \}$. 
When $|\Omega| > 1$, there are nondiagonal orbits too, so the rank is at least two, with equality if and only if $G$ acts transitively on $\Omega^{(2)}$.

For a nonnegative integer $m$, write 
\[
(m)_k = m(m-1) \cdots (m-k+1)
\]
for the falling factorial.
A permutation $g$ fixes $\fix(g)^k$ elements of $\Omega^k$ and $(\fix(g))_k$ elements of $\Omega^{(k)}$, so Lemma~\ref{lem:orbit-counting} applied to these actions shows that
\begin{gather*}
\text{the number of orbits of $G$ on $\Omega^k$} = \frac{1}{|G|} \sum_{g \in G} \fix(g)^k \\
\text{the number of orbits of $G$ on $\Omega^{(k)}$} = \frac{1}{|G|} \sum_{g \in G} (\fix(g))_k.
\end{gather*}

We say that the action is \emph{$k$-transitive} if the action on $\Omega^{(k)}$ is transitive, and \emph{sharply $k$-transitive} if the action on $\Omega^{(k)}$ is regular.
Said otherwise, the action is (sharply) $k$-transitive if given any two sequences of $k$ distinct elements of $\Omega$, say $(\alpha_1, \dots, \alpha_k)$ and $(\beta_1, \dots, \beta_k)$, there is a (unique) element $g \in G$ such that $(\alpha_1 g, \dots, \alpha_k g) = (\beta_1, \dots, \beta_k)$.
By the Orbit--Stabiliser Theorem, a $k$-transitive permutation group $G \leq \Sym(n)$ has order divisible by $(n)_k$, with equality if and only if it is sharply $k$-transitive.
Being sharply $1$-transitive is the same as being regular.

\begin{ex}
\label{ex:pgl2}
{\rm
Let $F$ be a field and let $d \geq 2$.
The group $\PGL_d(F)$ acts on the set of $1$-spaces of $F^d$.
When $d = 2$, this set is the \emph{projective line} over $F$, which we identify with $F \cup \{\infty\}$ (by mapping $\langle (x, 1) \rangle \mapsto x$ for $x \in F$ and $\langle (1, 0) \rangle \mapsto \infty$).
The group $\PGL_2(F)$ acts sharply $3$-transitively on the projective line (exercise).
}
\end{ex}

A permutation group $G \leq \Sym(\Omega)$ is a \emph{Frobenius group} if it is transitive and for all $\alpha \in \Omega$ the stabiliser $G_\alpha$ is nontrivial but for all distinct $\alpha, \beta \in \Omega$ the stabiliser $G_{\alpha \beta} = G_{\alpha} \cap G_{\beta}$ is trivial.
All sharply $2$-transitive groups of degree $n \geq 3$ are Frobenius.
Note that regular groups are excluded from being Frobenius groups, but a transitive group is \emph{regular or Frobenius} exactly when all its two-point stabilisers are trivial, that is, when no nontrivial element fixes more than one point, and this is the condition we will most often meet.
By the Orbit--Stabiliser Theorem, a regular or Frobenius group $G \leq \Sym(n)$ of rank $r$ has order $n(n-1)/(r-1)$.

\subsection{Coset actions}
\label{ss:permutations_equivalence}

Let $G$ be a group.
For a subgroup $H \leq G$, we write $G/H$ for the set of right cosets of $H$ in $G$.
The \emph{coset action} of $G$ on $G/H$ is as follows: for $g \in G$ and $Hx \in G/H$, we define $(Hx)g = H(xg)$.
This is a well-defined transitive action where the stabiliser of $Hx$ is $H^x$ and the kernel is $\bigcap_{x \in G} H^x$ (exercise).
In particular, $g \in G$ fixes a point of $G/H$ if and only if $g$ is \emph{conjugate into} $H$, that is, $g^x \in H$ for some $x \in G$.

Essentially, every transitive action of $G$ arises as a coset action. 
To make this precise, we need to define what we mean for two actions to be ``essentially the same''.

Two actions $G \to \Sym(\Omega_1)$ and $G \to \Sym(\Omega_2)$ are \emph{equivalent} if there is a bijection $\psi \colon \Omega_1 \to \Omega_2$ such that for all $g \in G$ and all $\omega \in \Omega_1$ we have
\[
(\omega g)\psi = (\omega \psi) g.
\]
Informally, two actions are equivalent if they are the same but with different labels for the points.
We will usually not distinguish between equivalent actions.

We now have the key to turning a question about group actions (or permutation groups) into one about abstract groups.

\begin{theo}
\label{thm:coset_action}
Let $G$ be a group acting transitively on a set $\Omega$. 
Fix $\omega \in \Omega$.
Then this action is equivalent to the coset action of $G$ on $G/G_\omega$.
\end{theo}

\begin{proof}
This is an improvement on the Orbit--Stabiliser Theorem. 
Define
\[
\psi\colon G/G_\omega \to \Omega \qquad G_\omega x \mapsto \omega x
\]
which we know is a bijection (noting that $\Omega = \omega G$ here, since the action is transitive). 
Moreover, for all $x, g \in G$ we have 
\[
(G_\omega xg)\psi = \omega (xg) = (\omega x)g = ((G_\omega x) \psi) g.
\]
Therefore, the two actions in the statement are equivalent.
\end{proof}

\subsection{Primitivity}
\label{ss:permutations_primitivity}

Let $G$ act on a set $\Omega$.
Rephrasing the definition, this action is \emph{transitive} if $G$ stabilises no proper nonempty subset of $\Omega$, and we say that $G$ is \emph{intransitive} otherwise.
Transitivity is our first ``irreducibility'' condition for actions.
Primitivity is our second.

A transitive action is \emph{primitive} if $|\Omega| > 1$ and $G$ stabilises no nontrivial partition of $\Omega$ (where a partition is trivial if there is one part or each part has size one) and \emph{imprimitive} otherwise.

As is standard, we only use the terms primitive and imprimitive in the context of transitive actions.
That said, the requirement that a primitive action be transitive is almost redundant: the only intransitive permutation group $G \leq \Sym(n)$ that stabilises no nontrivial partition is $1 \leq \Sym(2)$ (exercise).

As is \emph{not so standard}, we insist that $|\Omega| > 1$ for primitive actions, so $\Sym(1)$ is not a primitive permutation group.
This makes many theorems about primitive groups cleaner (e.g. Lemma~\ref{lem:primitive_maximal}), and it reflects the usual practice of excluding $1$ from being a prime number, the trivial group from being a simple group, the empty topological space from being connected, etc.

\begin{exer}
\label{exc:primitive_2-trans}
{\rm
Prove that every $2$-transitive group is primitive.
}
\end{exer}

We mention two key properties of primitive actions.

\begin{lem} 
\label{lem:primitive_normal}
Let $G \leq \Sym(\Omega)$ be transitive and let $N \trianglelefteqslant G$.
Then $G$ stabilises the partition of $\Omega$ given by the orbits of $N$.
In particular, any nontrivial normal subgroup of a primitive group is transitive.
\end{lem}

\begin{proof}
Let $\Pi = \{ \omega N \mid \omega \in \Omega \}$.
Let $g \in G$.
Since $N \trianglelefteqslant G$, we know that $Ng = gN$, so, for all $\omega \in \Omega$,
\[
(\omega N)g = \omega (Ng) = \omega (gN) = (\omega g) N.
\]
Therefore, $g$ stabilises $\Pi$.

Now assume that $G$ is primitive and $N$ is nontrivial.
Since $G$ is primitive and stabilises $\Pi$, the partition $\Pi$ must be trivial.
Since $G$ is transitive, the parts of $\Pi$ all have the same size, which is greater than one since $N$ is nontrivial, so, to be trivial, $\Pi$ must have exactly one part.
Said otherwise, $N$ has exactly one orbit, so $N$ is transitive.
\end{proof}

\begin{rem}
\label{rem:primitive_normal}
{\rm
The converse to Lemma~\ref{lem:primitive_normal} is not true. 
For example, any faithful transitive action of a simple group has the property that every nontrivial normal subgroup is transitive, but these need not be primitive.
(Indeed, by Lemma~\ref{lem:primitive_maximal} below, it suffices to consider the coset action of a nonabelian simple group $G$ on $G/H$ for a nonmaximal subgroup $H$ of $G$, e.g. $H = 1$.)
Transitive permutation groups where every nontrivial normal subgroup is transitive are called \emph{quasiprimitive}.
As Neumann points out in \cite{ref:Neumann06} (see also \cite[Section~I.8.6]{ref:Neumann11}), the two inequivalent properties of being primitive and quasiprimitive seem to be conflated by Galois in his \emph{Second M\'emoire}, written in 1830, both under the term \emph{primitif}.
By Jordan's \emph{Trait\'e} of 1870, the modern definition of primitivity had emerged with quasiprimitivity (not that it was called that then) being seen as a consequence. 
}
\end{rem}

\begin{lem} 
\label{lem:primitive_maximal}
Let $G$ act transitively on $\Omega$ and let $\omega \in \Omega$.
Then $G$ acts primitively on $\Omega$ if and only if $G_\omega$ is a maximal subgroup of $G$.
\end{lem}

\begin{proof}[Proof sketch]
First assume that $G$ acts imprimitively on $\Omega$.
Let $\Pi$ be a nontrivial partition of $\Omega$ stabilised by $G$ and let $\Delta$ be the part of $\Pi$ containing $\omega$.
Then the stabiliser $H$ of the set $\Delta$ satisfies $G_\omega < H < G$ (exercise).

Conversely, assume that there exists $G_\omega < H < G$.
Let $\Delta = \omega H$ and let $\Pi = \{ \Delta g \mid g \in G \}$.
Then $\Pi$ is a nontrivial partition stabilised by $G$ (exercise).
\end{proof}

\begin{rem}
\label{rem:primitive_maximal}
{\rm
By Theorem~\ref{thm:coset_action}, every transitive action of $G$ is equivalent to the coset action of $G$ on $G/H$ for a subgroup $H \leq G$. Lemma~\ref{lem:primitive_maximal} asserts that every primitive action of $G$ is equivalent to the coset action of $G$ on $G/H$ for a \emph{maximal} subgroup $H \leq G$.
}
\end{rem}

\begin{ex}
\label{ex:k-subsets}
{\rm
For $1 \leq k < n/2$, the actions of $S_n$ and $A_n$ on the $k$-subsets of $\{1, \dots, n\}$ are primitive, but for $n = 2k > 2$ they are not, since the partition into complementary pairs is invariant.
(As the degree is not $n$, we write $S_n$ and $A_n$, not $\Sym(n)$ and $\Alt(n)$.)
For example, for $A_4$ acting on the set of six $2$-subsets of $\{1, 2, 3, 4\}$, the invariant partition has three parts, which are the orbits of the normal subgroup $V_4 = \langle (1 \ 2)(3 \ 4), (1 \ 3)(2 \ 4) \rangle$, as in Lemma~\ref{lem:primitive_normal}.
}
\end{ex}

\subsection{Simple groups}
\label{ss:permutations_simple}

We now give a modern framework for studying finite primitive permutation groups, beginning with the special case of finite $2$-transitive groups (see Exercise~\ref{exc:primitive_2-trans}), for which we have the following theorem of Burnside \cite[Chapter~IX, Theorem~IX]{ref:Burnside97}.

\begin{theo}
\label{thm:burnside}
A finite $2$-transitive group is affine or almost simple.
\end{theo}

Let $V = \F_q^d$.
The \emph{affine general linear group} $\AGL_d(q) = V{:}\GL_d(q)$ acting on $V$ consists of all affine maps $v \mapsto u + vh$ for $u \in V$ and $h \in \GL_d(q)$.
In general, an \emph{affine group} is $V{:}G_0$ for $G_0 \leq \GL_d(q)$, which is primitive if and only if $G_0$ is irreducible (exercise).

\begin{exer}
\label{exc:agl1}
{\rm
Prove that $\AGL_1(q)$ is sharply $2$-transitive.
}
\end{exer}

A group $G$ is \emph{almost simple} if $T \leq G \leq \Aut(T)$ for a nonabelian simple group $T$.
For example, $A_d$ and $S_d$ are almost simple for $d \geq 5$, and so is $\PGL_2(q)$ for $q \geq 4$.
Exercise~\ref{exc:agl1} and Example~\ref{ex:pgl2} therefore realise the two cases of Theorem~\ref{thm:burnside}.

Burnside's formulation is that a $2$-transitive group $G \leq \Sym(\Omega)$ has a unique minimal normal subgroup $N$, which is either elementary abelian and regular, or nonabelian and simple.
In the first case, identify $\Omega$ with $N$, so $N$ acts by translations and $G_\omega$ by conjugation and $G = N{:}G_\omega$ is affine.
In the second case, $C_G(N) = 1$, by the uniqueness of $N$, so $G \leq \Aut(N)$.

For general primitive permutation groups, Burnside's formulation fails in two ways: there may be two minimal normal subgroups and a nonabelian minimal normal subgroup need not be simple.
Nonetheless, analysing the minimal normal subgroups still leads to a description of finite primitive permutation groups.
Here, the \emph{socle} of a group is the product of its minimal normal subgroups.

\begin{theo}[O'Nan--Scott Theorem]
\label{thm:onan_scott}
Let $G$ be a finite primitive permutation group.
Then one of the following holds
\begin{itemize}
\item[{\rm (a)}]
$G$ is affine
\item[{\rm (b)}]
$G$ is almost simple
\item[{\rm (c)}]
$G$ is of diagonal type
\item[{\rm (d)}]
$G$ is of product type, built from a group in (b) or (c)
\item[{\rm (e)}]
$G$ has a regular nonabelian socle.
\end{itemize}
\end{theo}

We will not need the precise definitions of the groups in (c)--(e), but let us explain how they arise.
The socle $N$ of $G$ is always $T^k$ for a simple group $T$ \cite[Theorem~4.3B]{ref:DixonMortimer96}, and the cases are distinguished by the projection $P$ of $N_\omega$ to a simple direct factor $T$ of $N$.
If $P = 1$, then $N$ is regular, giving (a) if $N$ is abelian and (e) if not, where the groups turn out to be \emph{twisted wreath products} (this uses the Classification of Finite Simple Groups, below).
If $P = T$, then $N_\omega$ is a diagonal subgroup of $N$, or a direct product of such, giving (c) or a product action of (c).
If $1 < P < T$, then we are in (b) if $k = 1$ or a product action of (b) if $k > 1$.

Theorem~\ref{thm:onan_scott} was announced independently by O'Nan and by Scott at the Santa Cruz conference in 1979, but only Scott's paper appeared in the proceedings \cite{ref:Scott80} (see \cite[p.~100]{ref:Cameron99}).
Case~(e) was originally omitted, and the first complete proof is due to Liebeck, Praeger and Saxl \cite{ref:LiebeckPraegerSaxl88}.

\begin{rem}
\label{rem:onan_scott}
{\rm
There are various formulations of the O'Nan--Scott Theorem. 
In one direction, Cameron \cite[Section~4.5]{ref:Cameron99} defines a third irreducibility condition for actions: an action is \emph{basic} if it preserves no Cartesian product structure.
He then gives a version of the O'Nan--Scott Theorem for basic permutation groups, which only features (a)--(c).
In another direction, Praeger \cite{ref:Praeger93} gives a version of the O'Nan--Scott Theorem with eight disjoint labelled cases that refine (c) and (d), and she extends this to a version that covers quasiprimitive groups.
}
\end{rem}

Theorems~\ref{thm:burnside} and~\ref{thm:onan_scott}, as stated, are elementary.
However, for most uses, one needs to know the actual groups that can arise. 
For this, we require the Classification of Finite Simple Groups (CFSG).

\begin{theo}[CFSG]
\label{thm:cfsg}
Let $G$ be a finite simple group.
Then one of the following holds
\begin{itemize}
\item[{\rm (a)}] 
$G$ is a cyclic group $C_p$ for a prime $p$
\item[{\rm (b)}]
$G$ is an alternating group $A_d$ for $d \geq 5$
\item[{\rm (c)}]
$G$ is a group of Lie type
\item[{\rm (d)}]
$G$ is one of $26$ sporadic simple groups.
\end{itemize}
\end{theo}

The prototypical example of a group of Lie type is $\PSL_d(q)$, and it is the only one we consider in any detail in this survey.
The smallest sporadic simple group is the Mathieu group ${\rm M}_{11}$, which is sharply $4$-transitive on $11$ points but also has a $3$-transitive action on $12$ points, with the latter playing a prominent role in Section~\ref{ss:derangements_prime}.

The proof of the CFSG was announced by Gorenstein in 1983, but it was only completed in 2004 with the treatment of the quasithin groups by Aschbacher and Smith \cite{ref:Aschbacher04}.
A second-generation proof by Gorenstein, Lyons and Solomon is still in progress, but in October 2026, as this survey was being finalised, a complete machine-checked formalisation of the first-generation proof in Lean was announced \cite{ref:Axiom26}.

As emphasised, for example, by Cameron \cite{ref:Cameron81}, the power of Theorems~\ref{thm:burnside} and~\ref{thm:onan_scott} comes from pairing them with the CFSG, which turns questions about $2$-transitive and primitive groups into questions about simple groups.
For instance, the $2$-transitive groups are now completely known \cite{ref:Cameron81, ref:DixonMortimer96}, and we use this classification in Section~\ref{ss:derangements_proportion_2}.

\clearpage

\section{Derangements}
\label{s:derangements}

\subsection{Origins}
\label{ss:derangements_history}

In his 1708 book ``Essay d’Analyse sur les Jeux de Hazard'' \cite{ref:Montmort08}, de~Montmort systematically studied games of chance that were popular at the time.
His discussion included a simplified version of the game \emph{jeu de treize}, which can be summarised as follows:

\begin{quote}
Pierre receives a suit of a standard deck of cards (with values 1, 2, \dots, 13) face down.
He deals the cards out one by one, calling out ``un'' as he turns over the first, ``deux'' as he turns over the second, etc.
Pierre loses if he never calls out the value of the card he turns over.
\end{quote}

De Montmort was interested in the probability that Pierre loses, which he asserted to be
\[
\frac{1}{0!} - \frac{1}{1!} + \frac{1}{2!} - \cdots - \frac{1}{13!} = \frac{63633137}{172972800} \approx \frac{1}{e}.
\]

Johann Bernoulli and Nikolaus Bernoulli contacted de~Montmort about the work in his book, and the 1713 edition of the book \cite{ref:Montmort13} features two proofs that the latter gave of the above assertion (in the general setting of $n$ cards).
One of these proofs is particularly notable, since it was the first application of the idea of Inclusion--Exclusion, which would be generalised soon after by de~Moivre (see \cite[Chapters~6 and~13]{ref:WilsonWatkins13}).
Tak\'acs \cite{ref:Takacs80} gives a good account of these developments.

The mathematical heart of this simplified version of \emph{jeu de treize} is the \emph{problem of coincidences} (or \emph{matches} or \emph{rencontres}).
This can be dressed up in any number of other ways.
For example, the \emph{hat check problem} asks: If everyone hands in their hat on the way in and receives a random hat on the way out, then what's the chance that nobody gets their own hat back? (see \cite{ref:Scoville66}).
However, behind all of these guises, there is one mathematical question, namely:

\begin{quote}
How many permutations fix no points?
\end{quote}

A \emph{derangement} is a permutation with no fixed points.
(The earliest use of this term with this meaning that the author could find is MacMahon's 1912 paper \cite{ref:MacMahon12}.)

The following is a modern statement of the theorem asserted by de~Montmort and proved by N.~Bernoulli.

\begin{theo}
\label{thm:de_montmort}
Let $n \geq 1$. 
The number of derangements of $n$ points is
\[
\sum_{i=0}^{n} (-1)^i \frac{n!}{i!},
\]
which is the closest integer to $n!/e$.
\end{theo}

\begin{proof}
The total number of permutations is $n!$.
Of them, exactly $(n-1)!$ fix a given point, and, more generally, exactly $(n-i)!$ fix $i$ given points.
Since the derangements are the permutations that fix no points, an Inclusion--Exclusion argument gives the number of derangements as
\[
n! - \binom{n}{1} \cdot (n-1)! + \binom{n}{2} (n-2)! - \cdots + (-1)^n \binom{n}{n} (n-n)!,
\]
or, more formally,
\[
\sum_{i=0}^{n} (-1)^i \binom{n}{i} (n-i)! = \sum_{i=0}^{n} (-1)^i \frac{n!}{i!},
\]
as required.
This is the closest number to $n!/e$ since
\begin{gather*}
\left|\frac{n!}{e} - \sum_{i=0}^{n} (-1)^i \frac{n!}{i!}\right| 
= \left|\sum_{i=0}^{\infty} (-1)^i \frac{n!}{i!} - \sum_{i=0}^{n} (-1)^i \frac{n!}{i!}\right| \\
= \left|\sum_{i=n+1}^{\infty} (-1)^i \frac{n!}{i!}\right| 
< \frac{n!}{(n+1)!} 
= \frac{1}{n+1} 
\leq \frac{1}{2}. \qedhere
\end{gather*}
\end{proof}

There are several different ways of counting derangements.
For example, \cite[Section~4.2]{ref:Cameron17} uses a recurrence, \cite[Section~8.1]{ref:Cameron17} uses the formal Principle of Inclusion--Exclusion and \cite[Exercise~10.6]{ref:Cameron17} adopts a beautiful approach via generating functions.

\subsection{Existence of derangements}
\label{ss:derangements_existence}

We now turn from the symmetric group and look at permutation groups in general.

In his 1872 paper ``Recherches sur les substitutions'' \cite{ref:Jordan72}, Jordan studied the sharply $k$-transitive groups with $k \geq 3$, giving a complete classification when $k \geq 4$ (the only examples are $\Alt(n)$ and $\Sym(n)$ on $n$ points and the Mathieu groups ${\rm M}_{11}$ and ${\rm M}_{12}$ on $11$ and $12$ points, respectively).
The first theorem in this paper is the following fundamental result. 
The proof we give here is not Jordan's original proof.

\begin{theo}
\label{thm:existence}
Every nontrivial finite transitive permutation group has a derangement.
\end{theo}

\begin{proof}
Fix $n > 1$ and a transitive permutation group $G \leq \Sym(n)$.
By the Orbit-Counting Lemma, 
\[
\frac{1}{|G|} \sum_{g \in G} \fix(g) = 1.
\]
The identity fixes $n > 1$ points, so there must be an element $g \in G$ with $\fix(g) < 1$, which must be a derangement. 
\end{proof}

\begin{rem}
\label{rem:existence_hypotheses}
{\rm
All of the hypotheses of Theorem~\ref{thm:existence} are necessary.
\begin{itemize}
\item[{\rm (a)}] 
\emph{Nontrivial.}
The trivial group $\Sym(1)$ is transitive on $\{1\}$ but has no derangements, since its only element fixes $1$.
\item[{\rm (b)}]
\emph{Finite.}
The finitary symmetric group $\FSym(\Nat)$, consisting of the permutations of $\Nat$ with finite support, is transitive on $\Nat$ but has no derangements, since every element fixes all but finitely many points.

As another example, for $d \geq 2$ the general linear group $\GL_d(\Com)$ acts transitively on the set $\Omega$ of $1$-spaces of $\Com^d$, but the induced permutation group $\PGL_d(\Com) \leq \Sym(\Omega)$ has no derangements since every matrix $A \in \GL_d(\Com)$ has an eigenvector $v$ (since the characteristic polynomial of $A$ has a root) and hence fixes the subspace $\langle v \rangle$.

Nevertheless, one can consider the class $\mathcal{X}$ of groups $G$ such that every transitive action of $G$ on a set of size greater than $1$ has a derangement.
Wiegold first studied this class in the 1970s \cite{ref:Wiegold76, ref:Wiegold77}, showing, for example, that it includes all virtually soluble groups (see Section~\ref{ss:connections_abstract}).
\item[{\rm (c)}]
\emph{Transitive.}
For $n \geq 2$, the stabiliser of $n$ in $\Sym(n)$ has two orbits but has no derangements, since every element fixes $n$.

Similarly, the affine group $\AGL_1(p)$ has two orbits on $\Omega = \F_p \sqcup \F_p^\times$, acting on $\F_p$ by affine maps and on $\F_p^\times$, with kernel $\F_p$, via the regular action of $\F_p^\times$, but has no derangements, since every translation fixes $\F_p^\times$ pointwise and every other element fixes a point of $\F_p$.

Anzanello and Eberhard \cite[Corollary~2]{ref:AnzanelloEberhard} proved that if $n \geq 3$ and $G \leq \Sym(n)$ has exactly two orbits and their lengths $a$ and $b$ satisfy $\gcd(a, b-1) = \gcd(a-1, b) = 1$, then $G$ has a derangement.
The case $a = b$ was conjectured by Ellis and Harper \cite{ref:EllisHarper26}, who proved it in certain cases (e.g. soluble groups, simple groups and groups of order at most $50000$).
\end{itemize}
}
\end{rem}

\subsection{Proportion of derangements: elementary bounds}
\label{ss:derangements_proportion_1}

The proof of Theorem~\ref{thm:existence} forces not just the existence of a single derangement but at least $n-1$ of them. 
This was, in fact, crucial for Jordan's application.
Can we say anything stronger? 

For a finite permutation group~$G$, write
\[
\delta(G) = \frac{|\{ g \in G \mid \text{$g$ is a derangement}\}|}{|G|}. 
\]
We begin with general bounds on $\delta(G)$.
The lower bound was proved by Cameron and Cohen in 1992 \cite{ref:CameronCohen92} and the upper bound by Guralnick, Isaacs and Spiga in 2015 \cite{ref:GuralnickIsaacsSpiga15}.
Given their short elementary proofs, it is remarkable that no such bounds in general transitive groups had been given in the 120 years between Jordan's paper in 1872 and Cameron and Cohen's paper in 1992. 
\nopagebreak

\begin{theo}
\label{thm:proportion}
Let $G \leq \Sym(n)$ be transitive of rank~$r$.
Then
\[
\frac{r-1}{n} \leq \delta(G) \leq \frac{r-1}{r}.
\]
The lower bound is met if and only if $G$ is regular or Frobenius, and the upper bound is met if and only if $G$ is regular.
\end{theo}

\begin{proof}
Write $\Omega = \{ 1, \dots, n \}$.
We begin with two applications of the Orbit-Counting Lemma.
Since $G$ has one orbit on~$\Omega$, 
\begin{equation}
\label{eq:sum-fix-1}
\sum_{g \in G} \fix(g) = |G|.
\end{equation}
Since $G$ has $r$ orbits on~$\Omega \times \Omega$ and $\fix(g)^2$ is the number of fixed points of~$g \in G$ on~$\Omega \times \Omega$, 
\begin{equation}
\label{eq:sum-fix-2}
\sum_{g \in G} \fix(g)^2 = r|G|.
\end{equation}

We now turn to the bounds.
Let $D$ be the set of derangements in~$G$ and let $C = G \setminus D$.

For the lower bound, consider the quantity
\[
f(g) = (\fix(g) - 1)(\fix(g) - n) = \fix(g)^2 - (n+1)\fix(g) + n.
\]
On the one hand, $f(g) \leq 0$ if $g \in C$, and $f(g) = n$ if $g \in D$, so
\[
\sum_{g \in G} f(g) \leq n|D|.
\]
On the other hand, applying \eqref{eq:sum-fix-1} and \eqref{eq:sum-fix-2} gives
\begin{gather*}
\sum_{g \in G} f(g) = \sum_{g \in G} \fix(g)^2 - (n+1) \sum_{g \in G} \fix(g) + \sum_{g \in G} n = (r - 1)|G|.
\end{gather*}
Therefore, $n|D| \geq (r-1)|G|$, that is,
\[
\delta(G) = \frac{|D|}{|G|} \geq \frac{r-1}{n}.
\]
Inspecting this argument, we get equality if and only if $f(g) = 0$ for all $g \in C$.
Since $(\fix(g) - n) = 0$ if and only if $g = 1$, we get equality if and only if $(\fix(g) - 1) = 0$ for all $g \in C \setminus \{1\}$, that is, if and only if $G$ is regular or Frobenius.

For the upper bound, we apply the Cauchy--Schwarz Inequality, to see that
\[
\biggl(\sum_{g \in C} \fix(g) \biggr)^2 \leq |C| \sum_{g \in C} \fix(g)^2.
\]
Since $\fix(g) = 0$ for all $g \in D$, by \eqref{eq:sum-fix-1}, 
\[
\biggl(\sum_{g \in C} \fix(g) \biggr)^2 = \biggl(\sum_{g \in G} \fix(g)\biggr)^2 = |G|^2,
\]
and, by \eqref{eq:sum-fix-2},
\[
\sum_{g \in C} \fix(g)^2 = \sum_{g \in G} \fix(g)^2 = r|G|.
\]
Therefore, $|G|^2 \leq |C|r|G|$, that is,
\[
\delta(G) = \frac{|D|}{|G|} = 1 - \frac{|C|}{|G|} \leq 1 - \frac{1}{r} = \frac{r - 1}{r}.
\]
We get equality in the Cauchy--Schwarz Inequality if and only if $\fix(g)$ is constant on~$C$, but $\fix(g) = n$ if and only if $g = 1$, so we get equality if and only if $C = \{1\}$, that is, if and only if $G$ is regular.
\end{proof}

\begin{rem}
\label{rem:frobenius_count}
{\rm
Let $G$ be a regular or Frobenius group of degree $n > 1$ and rank $r$, so $|G| = n(n-1)/a$ for $a = r - 1$ (see Section~\ref{ss:permutations_tuples}).
Then $\delta(G) = a/n$ by Theorem~\ref{thm:proportion}.
This can be seen directly.
The $n$ point stabilisers in $G$ pairwise intersect trivially, so $G$ has exactly $|G| - 1 - n(|G|/n - 1) = n - 1$ derangements.
}
\end{rem}

The following corollary of Theorem~\ref{thm:proportion} was proved by Cohen in 1991 (see \cite[p.~74]{ref:BuhlerLenstraPomerance93}).
His motivation is discussed in Section~\ref{ss:connections_galois}.

\begin{cor}
\label{cor:proportion}
Let $G \leq \Sym(n)$ be transitive with $n > 1$.
Then $\delta(G) \geq 1/n$.
\end{cor}


The Orbit-Counting Lemma also turns Inclusion--Exclusion into a formula, which recovers Theorem~\ref{thm:de_montmort} when $G = \Sym(n)$.

\begin{lem}
\label{lem:inclusion_exclusion}
Let $G \leq \Sym(\Omega)$ with $|\Omega| = n$, and let $o_i$ be the number of orbits of $G$ on $\Omega^{(i)}$ (so $o_0 = 1$). 
Then
\[
\delta(G) = \sum_{i=0}^{n} \frac{(-1)^i o_i}{i!},
\]
and truncating the sum after an odd (respectively even) $i$ gives a lower (respectively upper) bound for $\delta(G)$.
\end{lem}

\begin{proof}
For an integer $m \geq 0$, the sum $\sum_{i=0}^{m} (-1)^i \tbinom{m}{i}$ is $1$ if $m = 0$ and $0$ otherwise.
Hence, for $g \in G$, the sum $\sum_{i} (-1)^i \tbinom{\fix(g)}{i}$ is $1$ if $g$ is a derangement and $0$ otherwise.
Summing over $g \in G$ and using $\tbinom{\fix(g)}{i} = (\fix(g))_i / i!$ and the identity in Section~\ref{ss:permutations_tuples} gives the formula.
The bounds follow since the partial sums of $\sum_{i} (-1)^i \tbinom{m}{i}$ alternate about their limit.
\end{proof}

\begin{ex}
\label{ex:proportion_3-transitive}
{\rm
If $G \leq \Sym(n)$ is $3$-transitive, then $o_1 = o_2 = o_3 = 1$, so Lemma~\ref{lem:inclusion_exclusion} gives $\delta(G) \geq 1 - 1 + 1/2 - 1/6 = 1/3$.
}
\end{ex}

Part (a) of the next example shows that the upper bound in Theorem~\ref{thm:proportion} is optimal even when regular groups are excluded.

\begin{ex}
\label{ex:proportion_sharply_transitive}
{\rm
Let $G \leq \Sym(n)$ be sharply $k$-transitive.
Then $o_i = 1$ for $i \leq k$ and $o_i = (n)_i/(n)_k$ for $i > k$, so Lemma~\ref{lem:inclusion_exclusion}, together with the vanishing of ${\sum_{i=0}^{n} (-1)^i \tbinom{n}{i}}$, gives
\[
\delta(G) = \sum_{i=0}^{k-1} \frac{(-1)^i }{i!} \left( 1 - \frac{(n)_i}{(n)_k} \right).
\]
\begin{itemize}
\item[{\rm (a)}]
For a sharply $3$-transitive group of degree $n$ and rank $r = 2$, such as $\PGL_2(q)$ acting on the projective line over $\F_q$ (see Example~\ref{ex:pgl2}), this gives $\delta(G) = (1 - 1/n)/2$, which tends to $1/2 = (r-1)/r$ as $n \to \infty$.
\item[{\rm (b)}]
The sharply $4$- and $5$-transitive Mathieu groups of degree $11$ and $12$ satisfy $\delta({\rm M}_{11}) = 23/66$ and $\delta({\rm M}_{12}) = 107/288$.
\end{itemize}
}
\end{ex}

\subsection{Derangements of prime-power order}
\label{ss:derangements_prime_power}

In short, not only do derangements always exist in finite transitive permutation groups, but there are many of them. 
Can we always find derangements of a certain order?
In 1981, Fein, Kantor and Schacher proved the following theorem \cite{ref:FeinKantorSchacher81}, which addresses this.

\begin{theo}
\label{thm:fein_kantor_schacher}
Let $G \leq \Sym(n)$ be transitive with $n > 1$.
Then $G$ has a derangement of prime-power order.
\end{theo}

Section~\ref{ss:connections_brauer} discusses the original motivation for Theorem~\ref{thm:fein_kantor_schacher}.
Theorem~\ref{thm:existence}, giving the existence of a derangement, has a short elementary proof, but Theorem~\ref{thm:fein_kantor_schacher}, giving the existence of a derangement of prime-power order, uses the CFSG in a major way.

\begin{proof}[Proof sketch of Theorem~\ref{thm:fein_kantor_schacher}]
Suppose otherwise and fix a minimal counterexample~$G$.
More precisely, assume that the theorem fails for $G$ but holds for transitive permutation groups of degree less than~$n$ and degree~$n$ but order less than~$|G|$.
Write $\Omega = \{ 1, \dots, n \}$.

\vspace{4pt}

(1a) \emph{Reduction to primitive groups.}
Suppose there is a nontrivial $G$\babelhyphen{nobreak}invariant partition~$\Pi$ of~$\Omega$.
Since $G$ is transitive on~$\Omega$, the induced group $G^\Pi$ is transitive on~$\Pi$.
As $1 < |\Pi| < n$, by hypothesis, there exists a derangement $x \in G^\Pi$ of $p$-power order.
Let $g \in G$ be an element inducing $x \in G^\Pi$, and by raising $g$ to the $p'$-part of its order, we may assume that $g$ also has $p$-power order.
Since $G$ has no derangements of $p$-power order, $g$ must fix a point of $\Omega$, which forces $x$ to fix a point of $\Pi$, which is a contradiction.

\vspace{4pt}

(1b) \emph{Reduction to simple groups.}
Suppose there is a proper nontrivial normal subgroup~$N$ of~$G$.
Since $G$ is primitive, by Lemma~\ref{lem:primitive_normal}, $N$ is transitive.
As $1 < |N| < |G|$, there exists a derangement in $N$ of prime-power order, which is a contradiction.

\vspace{4pt}

(2) \emph{Invoke the CFSG.}
By (1), $G$ is simple.
By the CFSG (Theorem~\ref{thm:cfsg}), one of the following holds
\vspace{-7pt}

\begin{itemize}\setlength{\itemsep}{-4pt}
\item[{\rm (a)}] 
$G = C_p$ for $p$ prime
\item[{\rm (b)}]
$G = A_d$ for $d \geq 5$
\item[{\rm (c)}]
$G$ is a group of Lie type (e.g. $G = \PSL_d(q)$)
\item[{\rm (d)}]
$G$ is one of $26$ sporadic simple groups.
\end{itemize}

(3) \emph{Analysis of the simple groups.}
In summary, $G$ acts on $\Omega$ primitively, so, by Theorem~\ref{thm:coset_action} and Lemma~\ref{lem:primitive_maximal}, $\Omega = G/H$ for a maximal subgroup~$H$ of~$G$, $G$ has no derangements of prime-power order, so every element of $G$ of prime-power order is conjugate into $H$ (see Section~\ref{ss:permutations_equivalence}), and $G$ is simple, so, by the CFSG, $G$ appears in (a)--(d) above.

\vspace{4pt}

(3a) \emph{$G = C_p$ for $p$ prime.}
This case is easy to eliminate.
Any generator of $G$ has prime order but is not contained in $H = 1$.

\vspace{4pt}

(3b) \emph{$G = A_d$ for $d \geq 5$.}
This case uses 19th-century technology.

First assume that $H$ is intransitive, which, by maximality, forces $H = (S_k \times S_{d-k}) \cap A_d$ for some $1 \leq k < d/2$. 
For each prime divisor $p$ of $d$, write $d = p^ab$ and note that 
if $p$ is odd, then the product of $b$ cycles of length $p^a$ lies in $A_d$, so $H$ contains a conjugate of it and $p^a$ divides $k$.
If $p = 2$ and $a \geq 2$, then the product of $2b$ cycles of length $2^{a-1}$ lies in $A_d$, so $2^{a-1}$ divides $k$.
Therefore, $k$ is divisible by $d$ or $d/2$, which contradicts $k < d/2$.

Next assume that $H$ is imprimitive, stabilising a partition with $1 < k < d$ parts.
By Bertrand's Postulate (proved by Chebyshev \cite{ref:Chebyshev52}), we may fix a prime $d/2 < p < d$.
Then $H$ contains a $p$-cycle, but a $p$-cycle cannot nontrivially permute the parts (of which there are $k \leq d/2$), nor can it fix all parts (which have size $d/k \leq d/2$), which is a contradiction.

Now assume that $H$ is primitive.
The proper subgroup $H$ must contain a $3$-cycle, but this contradicts a theorem of Jordan \cite{ref:Jordan71}.

\vspace{4pt}

(3c) \emph{$G$ is a group of Lie type.}
This case uses 20th-century technology, but begins with some elementary techniques.
To simplify matters, let us assume that $G = \PSL_d(q)$ where $d \geq 3$ and $q$ is a power of a prime $p$.

A matrix with Jordan form $[J_2(1), J_1(1)^{d-2}]$ (and its image in $\PGL_d(q)$) is called a \emph{transvection}.
The transvections form a conjugacy class of elements of order $p$ in $G$, so $H$ must contain one, say $t$. 
Let $K$ be the normal subgroup of $H$ generated by the conjugacy class $t^H$.
By the Frattini argument $H = KN_H(P)$ where $P$ is a Sylow $p$-subgroup of $K$.
Being a nontrivial $p$-group, $P$ has a proper nonzero fixed space $U$, which is stabilised by $N_H(P)$.

Write $q = p^e$.
A prime $r$ is a \emph{primitive prime divisor} of $p^m-1$ if $r$ divides $p^m-1$ but not $p^i-1$ for all $1 \leq i < m$.
By Zsigmondy's Theorem \cite{ref:Zsigmondy82}, $p^{ed}-1$ has a primitive prime divisor when $d \geq 3$ unless $(q, d) \in \{ (2, 6), (4, 3) \}$; for exposition, we exclude these cases and fix a primitive prime divisor $r$ of $p^{ed}-1$.
Then $r$ divides $q^d - 1$ and hence $|G|$, but $r$ divides neither $|\GL_i(q)|$ for $i < d$ nor $|\GL_d(q_0)|$ for a proper subfield $\F_{q_0}$ of $\F_q$, so if $r$ divides the order of a subgroup $K$ of $G$, then $K$ must be irreducible.

Combining this, $H$ must contain an element of order $r$, so $r$ divides $|H|$.
Since $H = K N_H(P)$ and $r$ does not divide the order of the reducible subgroup $N_H(P)$, it divides $|K|$.

Therefore, $K$ corresponds to an irreducible subgroup of $\SL_d(q)$ that is generated by transvections.
In view of this, the possibilities for $K$ are given by a 1979 theorem of Kantor \cite[Theorem~II]{ref:Kantor79} (noting that $O_p(K) = 1$ since $K$ is irreducible).
As Kantor notes, this theorem exploits some ``very difficult classification theorems'' (such as Aschbacher's 1977 paper \cite{ref:Aschbacher77}), but it (like steps (3a) and (3b)) does not use the CFSG itself.
Considering the possibilities for $K$ given in \cite[Theorem~II]{ref:Kantor79} and making use of the fact that $r$ divides $|K|$, we quickly deduce that $K$ is one of the following:
\[
\begin{array}{cc}
\hline
K                & \text{conditions} \\
\hline
\PSp_d(q)        & \text{$d$ even} \\
{\rm O}^-_d(q)   & \text{$d$ even, $q$ even} \\
S_{d+1}, S_{d+2} & \text{$d$ even, $q = 2$} \\
\PSU_d(q^{1/2})  & \text{$d$ odd, $q$ square} \\
\hline
\end{array}
\]
Since $K$ is normal in $H$ and $H$ is maximal in $G$, we know that $H = N_G(K)$ which is a subgroup of $\Aut(K)$ (for example, since $K$ can be seen to be absolutely irreducible).
In all cases, $|H|$ is not divisible by a primitive prime divisor of $q^{d-1}-1$ but $|G|$ is, which is a contradiction.

\vspace{4pt}

(3d) \emph{$G$ is a sporadic group.}
This case uses 21st-century technology, simplifying the argument in \cite{ref:FeinKantorSchacher81}.
We now have complete lists of the conjugacy classes and maximal subgroups of all sporadic simple groups, accessible, for example, via the Character Table Library in GAP \cite{ref:CTblLib,ref:GAP}.
(The conjugacy classes are in the ATLAS \cite{ref:ATLAS}, but the maximal subgroups were only fully classified when the list for the Monster was completed in \cite{ref:DietrichLeePopiel25}.
Determining the maximal subgroups of the larger sporadic groups uses the CFSG (see \cite[Section~2.3]{ref:Wilson17}), so unlike the rest of step (3), our argument in (3d) uses the CFSG, but the argument in \cite{ref:FeinKantorSchacher81} avoids it.)

With this in place, we quickly see that if $|H|$ has the same distinct prime factors as $|G|$ then $(G,H)$ is one of the following, where we specify $H$ up to conjugacy in $\Aut(G)$:
\[
\begin{array}{ccccccccc}
\hline
G & {\rm M}_{11} & {\rm M}_{12} & {\rm M}_{12} & {\rm M}_{24} & {\rm HS}     & {\rm McL}    & {\rm Co}_{2} & {\rm Co}_{3} \\
H & \PSL_2(11)   & {\rm M}_{11} & \PSL_2(11)   & {\rm M}_{23} & {\rm M}_{22} & {\rm M}_{22} & {\rm M}_{23} & {\rm M}_{23} \\
\hline
\end{array}
\]
Moreover, in all but the first case $G$ has more conjugacy classes of elements of order $2$ or $3$ than $H$ does, which eliminates these possibilities for $(G, H)$.
Therefore, $G = {\rm M}_{11}$ and $H = \PSL_2(11)$, in which case every element of $G$ of prime order is conjugate into $H$ (so ${\rm M}_{11} \leq \Sym(12)$ has no derangements of prime order), but $G$ has elements of order $4$ and $8$ whereas $H$ has neither, and this final contradiction completes the proof.
\end{proof}

The proof of Theorem~\ref{thm:fein_kantor_schacher} demonstrates many common techniques in modern permutation group theory.
It also highlights a now ubiquitous three-step approach for proofs:
\begin{itemize}\setlength{\itemsep}{-4pt}
\item[{\rm (1)}]
Reduce to finite simple groups.
\item[{\rm (2)}]
Invoke the Classification of Finite Simple Groups.
\item[{\rm (3)}]
Analyse the finite simple groups.
\end{itemize}
In the original proof of Fein, Kantor and Schacher, steps (1) and (3) do not use the CFSG, though it would be misleading to describe the ingredients in (3c) as elementary.
However, there are other proofs where (1) or (3) themselves use the CFSG (see the discussion of the proofs of Theorems~\ref{thm:malle_navarro_olsson} and~\ref{thm:giudici}, respectively).
More generally, the structure theorem in step (2) need not be the CFSG.
Indeed, in the proofs of Theorems~\ref{thm:guralnick_wan} and~\ref{thm:garzoni} below, it is Theorem~\ref{thm:burnside} on $2$-transitive groups or the O'Nan--Scott Theorem, respectively, and the CFSG enters only in the analysis in step (3).

\begin{rem}
\label{rem:p-elements}
{\rm
The lengths of all the cycles of a derangement of $p$-power order are divisible by $p$, so a transitive group of degree $n$ can have a derangement of $p$-power order only if $p \mid n$.
Conversely, if all cycles of $g$ have length divisible by $p$, then raising $g$ to the $p'$-part of its order gives a derangement of $p$-power order.
}
\end{rem}

By Theorem~\ref{thm:fein_kantor_schacher}, every transitive group of degree $n > 1$ has a derangement of $p$-power order for \emph{some} prime $p$ dividing $n$, but the next example shows that it need not have one for \emph{every} such prime.

\begin{ex}
\label{ex:isbell}
{\rm
Let $q$ be a prime power satisfying $q \equiv 1 \pmod{4}$.
Consider the action of $G = \PSL_2(q) \leq \PGL_2(q)$ on the projective line over $\F_q$.
Derangements in $G$ correspond to matrices in $\SL_2(q)$ that act irreducibly, and hence they have order dividing $(q+1)/2$ which is odd, so $G$ has no derangements of $2$-power order.
}
\end{ex}

If $n$ is a power of a prime, then the prime $p$ above is uniquely determined.
What about numbers that are close to being a power of a prime?
This question leads to the following conjecture, which was first made by Isbell in 1960 \cite{ref:Isbell60} for $p = 2$ and generalised to its current form by Cameron, Frankl and Kantor in 1989 \cite{ref:CameronFranklKantor89}.
Isbell's original motivation for this conjecture is discussed in Section~\ref{ss:connections_isbell}.

\begin{conj}[Isbell's Conjecture]
\label{conj:isbell}
For each prime $p$, there is a map $m_p \colon \Nat\to \Nat$ such that every transitive subgroup of $\Sym(p^a b)$ with $p \nmid b$ and $a \geq m_p(b)$ has a derangement of $p$-power order.
\end{conj}

This conjecture remains open.
However, Halpenny \cite{ref:Halpenny91} proved it for primitive groups, and Bereczky \cite{ref:Bereczky95} handled the case where $b < 3(p+1)/2$, taking $m_p(b) = 2$ if $b = p+1$ and $m_p(b) = 1$ otherwise.
Spiga \cite{ref:Spiga04} also established some lower bounds on the map $m_p$.
For example, if $b = (p^k-1)/(p-1)$, then $m_p(b) \geq k$ (for $k = 2$, compare with Bereczky's bound above and Example~\ref{ex:elusive} below), and see also \cite{ref:CameronEllisRaynaud19} for more in this direction.
As an approach to proving Isbell's Conjecture, Cameron suggested a stronger version, but this has since been refuted for $p \geq 3$ by Crestani and Spiga \cite{ref:CrestaniSpiga10} and Spiga \cite{ref:Spiga13}.

\subsection{Derangements of prime order}
\label{ss:derangements_prime}

Let $G \leq \Sym(n)$ be transitive with $n > 1$.
Then $G$ has a derangement of prime-power order, by Theorem~\ref{thm:fein_kantor_schacher},  but as the group ${\rm M}_{11} \leq \Sym(12)$ shows, $G$ need not have a derangement of prime order (see (3d) in the proof of Theorem~\ref{thm:fein_kantor_schacher}).

Following Cameron et~al.\ \cite{ref:CameronGiudiciJonesKantorKlinMarusicNowitz02}, we say that a transitive group $G \leq \Sym(n)$ with $n > 1$ is \emph{elusive} if it has no derangement of prime order.
As they explain
\begin{quote}
The name is intended to suggest that such groups are not easy to find. 
It also suggests that, given such a group $G$, we lack one of the standard tools for studying the $G$-invariant graphs [\dots]
\end{quote}

The second point is discussed in Section~\ref{ss:connections_polycirculant}.
On the first point, elusive groups do appear to be elusive!
This suggests the following problem posed by Cameron et~al.\ \cite[Problem~6.2]{ref:CameronGiudiciJonesKantorKlinMarusicNowitz02}.

\begin{que}
\label{que:elusive_density}
Does the set 
\[
\{ n \in \Nat \mid \text{there exists an elusive group of degree $n$} \}
\] 
have density zero in $\Nat$?
\end{que}

The following example is given in \cite{ref:FeinKantorSchacher81} (see also \cite{ref:CameronGiudiciJonesKantorKlinMarusicNowitz02}).

\begin{ex}
\label{ex:elusive}
{\rm
Let $p$ be a Mersenne prime, let $G$ be the affine group $\AGL_1(p^2)$, let $H$ be the subgroup $\AGL_1(p)$ and consider the coset action of $G$ on $G/H$ of degree $n = p(p+1)$.
By Remark~\ref{rem:p-elements}, it suffices to show that every element of $G$ of prime order dividing $n$ is conjugate into $H$, and since $p$ is a Mersenne prime, these primes are $2$ and $p$, both of which divide $|H|$.
The Sylow $2$-subgroups of $G$ are cyclic, so all elements of order $2$ are conjugate, and the Sylow $p$-subgroup is the normal subgroup $C_p^2$, on whose nontrivial elements $G$ acts transitively, so all elements of order $p$ are conjugate.
Hence, every element of order $2$ or $p$ is conjugate into $H$, so $G$ has no derangements of prime order on $G/H$.
As all elements of $p$-power order have order $p$, the group $G$ also has no derangements of $p$-power order.
}
\end{ex}

The paper \cite{ref:CameronGiudiciJonesKantorKlinMarusicNowitz02} presents many other constructions and more recent ones are surveyed in \cite[Sections~2 and~3]{ref:ArezoomandAbdollahiSpiga19}. 
Many of these have a similar flavour to the one we have just seen, with only the 2026 paper \cite{ref:ChenLeeMitrovicOBrienXia26} not relying on semidirect product constructions.
This paper gave the first example of elusive groups of odd degree, and the following is the smallest example in their infinite family. 

\begin{ex}
\label{ex:elusive_odd}
{\rm
The example is built from the degree~$15$ action of $A_5$ on $A_5/K$, where $K = \langle (1 \ 2)(3 \ 4), (1 \ 3)(2 \ 4) \rangle$.
Let $U = \F_3^4$ and $V = \F_5^3$ be the fully deleted permutation modules for $A_5$, regarded as the vectors in $\F_3^5$ and $\F_5^5$ with coordinate sum zero, the latter modulo constant vectors.
In \cite[Section~4]{ref:ChenLeeMitrovicOBrienXia26}, a nonsplit extension $G = (U \times V).A_5$ is constructed in which every element of order $3$ or $5$ lies in $U \times V$.
Let $U_0$ and $V_0$ be the hyperplanes of $U$ and $V$ defined by $x_1 + x_4 = x_2 + x_3$, which are $K$-invariant, and let $H = (U_0 \times V_0){:}K$.
Consider the degree~$225$ action of $G$ on $\Omega = G/H$.
As $225 = 3^25^2$ and every element of order $3$ or $5$ lies in $U \times V$, on which $G$ acts via $A_5$, to show that $G$ is elusive, it suffices to show that the $A_5$-translates of $U_0$ and $V_0$ cover $U$ and $V$, respectively.
These translates are the hyperplanes defined by $x_i + x_l = x_j + x_k$, where $\{i, l\}$ and $\{j, k\}$ are disjoint pairs, so we need every $x \in \F_p^5$ with coordinate sum~$0$, for $p \in \{3, 5\}$, to have two disjoint pairs of coordinates with equal sums (exercise).
}
\end{ex}

Returning to the proof of Theorem~\ref{thm:fein_kantor_schacher}, in most cases our analysis of the simple groups in step (3) actually produced a derangement of prime order, and the reduction in (1b) extends this to all primitive groups.
However, the reduction in (1a) does not extend it to all transitive groups: if $x \in G^\Pi$ is a derangement of prime order $p$, then we can lift $x$ to an element of $G$ of $p$-power order, but not necessarily to one of order $p$.
This failure to proceed inductively via quotients is what makes elusive groups so difficult to study.
Nevertheless, the argument does extend to primitive, and quasiprimitive, groups: in 2003, Giudici \cite{ref:Giudici03} proved the following.

\begin{theo}
\label{thm:giudici}
Let $G \leq \Sym(n)$ be elusive with a transitive minimal normal subgroup.
Then $G = {\rm M}_{11}$ and $n = 12$, or $G$ is of product type built from ${\rm M}_{11}$ and $n = 12^k$ for some $k \geq 2$.
\end{theo}

Step (3) of our earlier analysis already handles most cases.
In (3a), cyclic groups pose no problem, and in (3d), the only elusive primitive sporadic group was ${\rm M}_{11} \leq \Sym(12)$.
In (3b), the only action of $A_d$ in which we did not find a derangement of prime order was the one on $k$-subsets of $\{1, \dots, d\}$ for $1 \leq k < d/2$, but one exists there too (exercise).
In (3c), we also found a derangement of prime order, but we did not consider all groups of Lie type (nor are all cases discussed in \cite{ref:FeinKantorSchacher81}).
However, Giudici proceeds differently, and uses a classification by Liebeck, Praeger and Saxl \cite{ref:LiebeckPraegerSaxl00} of the maximal subgroups $H$ of an almost simple group $G$ with socle $T$ such that every prime divisor of $|G|$ divides $|H \cap T|$.
This classification relies on the CFSG, so, unlike in the original proof of Theorem~\ref{thm:fein_kantor_schacher}, step (3) uses the CFSG too.
In fact, Giudici's argument handles all almost simple groups.

\begin{theo}
\label{thm:giudici_almost_simple}
Let $G \leq \Sym(n)$ be elusive and almost simple.
Then $G = {\rm M}_{11}$ or $G = {\rm M}_{10}$, and $n = 12$.
\end{theo}

Any transitive subgroup of an elusive group is elusive, so the example ${\rm M}_{11} \leq \Sym(12)$ yields ${\rm M}_{10} = A_6.2$ and $\AGL_1(9)$ (the smallest case of Example~\ref{ex:elusive}).

We conclude with a tempting question, which no known example of an elusive group has yet refuted.

\begin{que}
\label{que:prime-square}
Let $G \leq \Sym(n)$ be transitive with $n > 1$.
Does $G$ have a derangement of prime or prime-squared order?
\end{que}

\subsection{Proportion of derangements: advanced bounds}
\label{ss:derangements_proportion_2}

With some experience of proofs using the CFSG, we now return to proportion and discuss some more advanced bounds.

Shortly after Cameron and Cohen's work, a paper of Boston et~al.\ \cite{ref:Boston93} made many interesting observations about the proportion of derangements, including the fact that the set $\{ \delta(G) \mid \text{$G$ is transitive} \}$ is dense in~$[0,1]$.
Their paper concludes with several provocative questions \cite[Section~7]{ref:Boston93}, the fourth of which focusses on refining the lower bound in Theorem~\ref{thm:proportion}, which we discuss now.

In Theorem~\ref{thm:proportion}, as $r \geq 2$, the lower bound implies the following.

\begin{cor}
\label{cor:proportion_frobenius}
Let $G \leq \Sym(n)$ be transitive with $n \geq 3$.
Then one of the following holds
\begin{itemize}
\item [{\rm (a)}]
$G$ is Frobenius of order $n(n-1)$ and $\delta(G) = 1/n$
\item [{\rm (b)}]
$\delta(G) > 1/n$.
\end{itemize}
\end{cor}

We exclude $n = 2$ since $\Sym(2)$ is regular, so, by definition, not a Frobenius group.

In 1997, Guralnick and Wan \cite[Theorem~1.3]{ref:GuralnickWan97} used the CFSG to strengthen Corollary~\ref{cor:proportion_frobenius} as follows.

\begin{theo}
\label{thm:guralnick_wan}
Let $G \leq \Sym(n)$ be transitive with $n \geq 10$.
Then one of the following holds
\begin{itemize}
\item[{\rm (a)}] 
$G$ is Frobenius of order $n(n-1)/a$ and $\delta(G) = a/n \leq 2/n$
\item[{\rm (b)}] 
$\delta(G) > 2/n$.
\end{itemize}
\end{theo}

Exceptions arise for $n \leq 9$, and Guralnick and Wan classified these.
(Two exceptions are missing from their list: $G = \PSL_3(2) \leq \Sym(7)$ and $G = \AGaL_1(9) \leq \Sym(9)$, for which $\delta(G) = 2/n$.)

The motivation for, and application of, Guralnick and Wan's theorem is discussed in Section~\ref{ss:connections_values}.

We now turn to the proof.

\begin{proof}[Proof sketch of Theorem~\ref{thm:guralnick_wan}] \quad

(1a) \emph{Reduction to $2$-transitive groups.} 
Let $r$ be the rank of $G$.
If $G$ is regular or Frobenius, then $\delta(G) = a/n$ with $a = r - 1$ by Remark~\ref{rem:frobenius_count}, which exceeds $2/n$ unless $a \leq 2$, in which case $G$ is Frobenius and (a) holds.
If $G$ is not regular or Frobenius, then Theorem~\ref{thm:proportion} gives $\delta(G) > (r-1)/n$, so $\delta(G) > 2/n$ unless $r = 2$, that is, unless $G$ is $2$-transitive.

\vspace{4pt}

(1b) \emph{Reduction to a bound on minimal degree.}
Let $t$ be the number of orbits of a two-point stabiliser.
Then
\[
\sum_{g \in G} \fix(g) = |G|, \quad \sum_{g \in G} \fix(g)^2 = 2|G|, \quad \sum_{g \in G} \fix(g)^3 = (t+2)|G|.
\] 
Let $\mu$ be the \emph{minimal degree} of $G$, the smallest number of points moved by a nontrivial element.
Consider $\sum_{g \in G} f'(g)$ for
\[
f'(g) = (\fix(g) - 1)(\fix(g) - n)(n - \mu - \fix(g)),
\]
noting that $f'(g) \leq 0$ if $g$ has a fixed point and $f'(g) = n(n-\mu)$ otherwise.
Arguing as in the proof of Theorem~\ref{thm:proportion} we deduce that
\[
\delta(G) \geq \frac{1}{n} + \left( \frac{n-t}{n-\mu} \right) \cdot \frac{1}{n},
\]
so it suffices to prove that $t < \mu$.

\vspace{4pt}

(2) \emph{Invoke Theorem~\ref{thm:burnside}.}
By (1), $G$ is $2$-transitive.
By Theorem~\ref{thm:burnside}: (a) $G$ is affine, or (b) $G$ is almost simple.

\vspace{4pt}

(3) \emph{Analysis of the $2$-transitive groups.}
In summary, $G$ is $2$-transitive of degree $n \geq 10$ and it suffices to show that $t < \mu$.

\vspace{4pt}

(3a) \emph{$G$ is affine.}
This case needs no classification.
Let $G = V{:}G_0$ with $V = \F_q^d$ and $G_0 \leq \GaL_d(q)$, and assume that the affine action of $G$ on $V$ is $2$-transitive.
We directly prove that $t < \mu$ when $n \geq 10$.

For example, let $d = 1$, $q = p^2$ for an odd prime $p$, and $G_0 = \GaL_1(p^2)$.
Let $g \in G$ be nontrivial with a fixed point.
Replacing $g$ by a conjugate, we may assume that $g \in G_0$.
If $g \in \GL_1(p^2)$, then $0$ is its unique fixed point. 
Otherwise, $g$ is the map $x \mapsto \lambda x^p$ for some $\lambda \in \F_{p^2}^\times$, whose fixed points are $0$ and the solutions of $x^{p-1} = \lambda^{-1}$, so it fixes at most $p$ points, with equality for the map $\phi \colon x \mapsto x^p$, whose fixed points are the elements of $\F_p$.
Therefore, $n - \mu = p$.
The orbits of $G_{01} = \langle \phi \rangle$ are the $p$ points of $\F_p$ and the $(p^2-p)/2$ pairs $\{x, x^p\}$ with $x \notin \F_p$, so $t = (p^2+p)/2$.
Therefore, if $p \geq 5$, then $t < p^2 - p = \mu$, while $G$ has exactly $(p+1)(p^2-1)$ derangements (exercise), so $\delta(G) = (p+1)/(2p^2) = (\sqrt{n}+1)/(2n)$, which is $2/9$ for $p = 3$.
Hence, $n \geq 10$ is best possible in this case.

\vspace{4pt}

(3b) \emph{$G$ is almost simple.}
This case uses the classification of $2$-transitive almost simple groups, and hence the CFSG.
We can leave aside $3$-transitive groups (see Example~\ref{ex:proportion_3-transitive}).
For the remaining groups, we can directly compute $\mu$ and $t$ and conclude that $t < \mu$ for all such groups with $n \geq 8$.

For example, let $d \geq 3$ and consider the $2$-transitive action of $G = \GL_d(2)$ on the set $\Omega$ of $n = 2^d-1$ nonzero vectors of $\F_2^d$.
The fixed points of a matrix are its $1$-eigenvectors, so the number of fixed points of a nonidentity matrix is maximised when the $1$-eigenspace is $(d-1)$-dimensional.
Hence $n - \mu = 2^{d-1}-1$.
The orbits of $G_{uv}$, for distinct nonzero $u, v \in \F_2^d$, are $\{ u \}, \{ v \}, \{u + v\}, \Omega \setminus \{ u, v, u + v \}$, so $t = 4$.
Therefore, if $d \geq 4$, then $t = 4 < 2^d - 2^{d-1} = \mu$, while for $\GL_3(2) = \PSL_3(2)$, acting on the $7$ nonzero vectors of $\F_2^3$, the derangements are the $48$ elements of order $7$, so $\delta(G) = 48/168 = 2/7$.
Hence, $n \geq 8$ is best possible in this case.
\end{proof}

Corollary~\ref{cor:proportion_frobenius} and Theorem~\ref{thm:guralnick_wan} suggest a natural pattern, and in 2025, Garzoni \cite{ref:Garzoni25} proved the following.

\begin{theo}
\label{thm:garzoni}
Let $G \leq \Sym(n)$ be transitive with $n$ sufficiently large.
Then one of the following holds
\begin{itemize}
\item[{\rm (a)}] 
$G$ is primitive Frobenius of order $n(n-1)/a$ and $\delta(G) = a/n$
\item[{\rm (b)}] 
$\delta(G) \geq (\sqrt{n}+1)/(2n)$.
\end{itemize}
\end{theo}

Theorem~\ref{thm:garzoni} confirms a conjecture of Guralnick and Tiep \cite{ref:GuralnickTiep03} (with a slight correction) for sufficiently large~$n$, and the bound in (b) is attained by $\AGaL_1(p^2)$ (see (3a) in the proof of Theorem~\ref{thm:guralnick_wan}).

Before discussing the proof of Theorem~\ref{thm:garzoni}, we first return to another question that was raised by Boston et~al.\ (the second question in \cite[Section~7]{ref:Boston93}).
This question asks for the minimum value of $\delta(G)$ as $G$ ranges over simple transitive groups.
Based on this question and an independent unpublished conjecture of Shalev (attribution to Shalev is given in \cite{ref:FulmanGuralnick04}), the Boston--Shalev Conjecture came to refer to the assertion that $\delta(G)$ is bounded away from zero as $G$ ranges over simple transitive groups.
This conjecture was proved by Fulman and Guralnick in a series of papers between~2003 and~2018 \cite{ref:FulmanGuralnick03,ref:FulmanGuralnick12,ref:FulmanGuralnick17,ref:FulmanGuralnick18}.
This proof uses the CFSG and detailed information about the maximal subgroups and conjugacy classes in simple groups, and is therefore beyond the scope of this survey.

\begin{theo}
\label{thm:boston_shalev}
There exists a constant $\epsilon > 0$ such that all finite simple transitive permutation groups $G$ satisfy $\delta(G) \geq \epsilon$.
\end{theo}

An explicit constant $\epsilon$ is not given by Fulman and Guralnick, but they do prove that $\epsilon$ can be taken to be $0.016$ if an unspecified finite collection of simple groups~$G$ are excluded.
Originally, Boston et~al.\ suggested that $2/7 \approx 0.286$ might be the optimal constant (which, as noted above, is achieved by the degree $7$ action of $\PSL_3(2)$), but Burness subsequently highlighted that the degree~$2925$ action of ${}^2F_4(2)'$ achieves $89/325 \approx 0.274$ and speculated that this might be the optimal bound (see, for example, \cite[Remark~1.3(c)]{ref:BurnessFusari25}).

\begin{que}
\label{que:boston_shalev}
Do all finite simple transitive permutation groups~$G$ satisfy 
\[
\delta(G) \geq \frac{89}{325}?
\]
\end{que}

What is the smallest value of $\delta(G)$ asymptotically as $G$ ranges over simple transitive groups?
Neumann and Praeger \cite[Theorem~4.3]{ref:NeumannPraeger98} proved that the probability that a matrix in $\GL_d(q)$ has no eigenvector tends to $\prod_{i = 1}^{\infty}(1-q^{-i})^{q-1}$ as $d \to \infty$.
This probability is the same as the proportion of elements of $\PGL_d(q)$ that are derangements in its action on the $1$-spaces of $\F_q^d$.
In particular, taking $q = 2$, for which $\PSL_d(2) = \PGL_d(2)$, gives $\delta(\PSL_d(2)) \to \prod_{i = 1}^{\infty}(1-2^{-i}) \approx 0.288$ as $d \to \infty$.
In 2016, Burness \cite{ref:Burness16} conjectured that this is the worst case.

\begin{conj}
\label{conj:boston_shalev_asymptotic}
Let $(G_n)$ be a sequence of finite simple transitive permutation groups such that $|G_n| \to \infty$ as $n \to \infty$. 
Then 
\[
\liminf_{n \to \infty} \delta(G_n) \geq \prod_{i = 1}^{\infty}(1- 2^{-i}).
\]
\end{conj}

We can now return to the proof of Garzoni's theorem.

\begin{proof}[Proof sketch of Theorem~\ref{thm:garzoni}]
Write $\Omega = \{ 1, \dots, n \}$ and, for ease of discussion, $f(n) = (\sqrt{n} + 1)/(2n)$.

\vspace{4pt}

(1) \emph{Reduction to primitive groups.}
Assume the result for primitive groups.
Let $G$ be imprimitive and let $\Pi$ be a nontrivial partition of $\Omega$ stabilised by $G$ with as few parts as possible, so $G^\Pi$ is primitive.
Write $k = |\Pi|$, so $n = mk$ where $m$ is the size of a part.
Every element of $G$ inducing a derangement in $G^\Pi$ is a derangement in $G$, so $\delta(G) \geq \delta(G^\Pi) \geq 1/k$, which suffices if $k$ is bounded, and if $G^\Pi$ is not Frobenius and $k$ is large, then $\delta(G) \geq \delta(G^\Pi) \geq f(k) > f(n)$.

Now assume that $G^\Pi$ is Frobenius, with point stabilisers of order $s \geq 2$, so $|G^\Pi| = ks$.
We count two sorts of derangements in $G$: those that fix no part and those that fix exactly one.
First, by Remark~\ref{rem:frobenius_count}, $G^\Pi$ has exactly $k-1$ derangements, so the elements of $G$ inducing one of them form a proportion $(k-1)/(ks)$ of $G$.
Second, if $H$ is the stabiliser in $G$ of a part $\Delta$ and $K$ is the kernel of the action on $\Pi$, then at least $|H|/m$ elements of $H$ are derangements on $\Delta$, by Corollary~\ref{cor:proportion}, and those outside $K$ fix no other part, so are derangements on $\Omega$. 
As $|H| = |G|/k$ and $|K| = |G|/(ks)$, counting over the $k$ parts gives a proportion at least $1/m - 1/s$.
Therefore, $\delta(G) \geq (k-1)/(ks) + 1/m - 1/s \geq 1/m - 1/(2k)$, which, together with $\delta(G) \geq 1/k$, gives $\delta(G) \geq \sqrt{2/(3n)} > f(n)$, as the larger of these bounds is smallest when $k = \sqrt{3n/2}$.

\vspace{4pt}

(2) \emph{Invoke the O'Nan--Scott Theorem.}
By (1), $G$ is primitive.
By the O'Nan--Scott Theorem (Theorem~\ref{thm:onan_scott}), we reduce to two cases: (a) $G$ is not affine, in which case its socle is $T^k$ for a nonabelian simple group $T$, or (b) $G$ is affine.

\vspace{4pt}

(3) \emph{Analysis of the primitive groups.}
In summary, $G$ is primitive of sufficiently large degree $n$ and we must show that $\delta(G) \geq f(n)$ unless $G$ is Frobenius with $\delta(G) < f(n)$.

\vspace{4pt}

(3a) \emph{$G$ is not affine.}
It suffices to show that $\delta(G) \geq c/\log n$ for an absolute constant $c > 0$.
We focus on almost simple groups $G$ with socle $T$, but similar arguments handle the other cases \cite[Section~8]{ref:FulmanGuralnick03}.
Here, $T$ is transitive, so $\delta(G) \geq \delta(T)/|G:T| \geq \epsilon/|\mathrm{Out}(T)| \geq \epsilon/(C\log n)$ for constants $\epsilon, C > 0$, by Theorem~\ref{thm:boston_shalev} and the bound $|\mathrm{Out}(T)| \leq C\log n$, which follows from the CFSG.

\vspace{4pt}

(3b) \emph{$G$ is affine.}
Let $G = V{:}G_0$ with $V = \F_q^d$.
The map $v \mapsto u + vh$ fixes a point if and only if $u$ lies in the image of $1 - h$, whose kernel is $C_V(h)$, so the proportion of derangements in the coset $Vh$ is $1 - 1/|C_V(h)|$, and we need many elements of $G_0$ with eigenvalue~$1$, which Garzoni obtains in a long argument exploiting the structure of irreducible linear groups \cite[Theorem~1.6]{ref:Garzoni25}.
\end{proof}

\clearpage

\section{Connections}
\label{s:connections}

\vspace{-11pt}

\subsection{Group theory: normal coverings \& invariable generation}
\label{ss:connections_abstract}

\sectionquestion{Can the conjugates of a proper subgroup cover a group?}

Theorem~\ref{thm:existence} will imply that the answer for finite groups is no.
For this, we need a group theoretic interpretation of derangements.

\begin{lem}
\label{lem:abstract}
Let $G$ be a group. 
Then the following are equivalent
\begin{itemize}
\item [{\rm (a)}]
every transitive action of $G$ of degree $n > 1$ has a derangement
\item [{\rm (b)}]
$G$ is not the union of the conjugates of any proper subgroup $H$
\item [{\rm (c)}]
every set with an element from each conjugacy class generates $G$.
\end{itemize}
\end{lem}

\begin{proof}
\textsc{not} (a) $\Leftrightarrow$ \textsc{not} (b).
If no element of $G$ is a derangement in a transitive action on $\Omega$ and $|\Omega| > 1$, then every element of $G$ is conjugate into $G_\omega$ for some $\omega \in \Omega$, that is, $G$ is the union of the conjugates of $G_\omega$, and $G_\omega < G$.
Conversely, if $G$ is the union of the conjugates of a proper subgroup $H$, then no element of $G$ is a derangement in the coset action on $G/H$ and $|G/H| > 1$.

\textsc{not} (b) $\Leftrightarrow$ \textsc{not} (c).
If $G$ is the union of conjugates of $H < G$, then for all $x \in G$ there exists $g_x \in G$ such that $x^{g_x} \in H$, so $X = \{ x^{g_x} \mid x \in G \}$ meets each conjugacy class but $\langle X \rangle \leq H < G$.
Conversely, if $X$ contains an element from each conjugacy class and $\langle X \rangle < G$, then every element of $G$ is conjugate to an element of $\langle X \rangle$, so $G$ is the union of conjugates of $\langle X \rangle$.
\end{proof}

We now discuss two corollaries of Theorem~\ref{thm:existence}.

\begin{cor}
\label{cor:covering}
No finite group is the union of conjugates of a proper subgroup.
\end{cor}

This leads to the following definition.
A \emph{normal covering} of a group $G$ is a set $\mathcal{H}$ of proper subgroups of $G$ such that
\[
G = \bigcup_{H \in \mathcal{H}}\bigcup_{g \in G} H^g,
\]
and the \emph{normal covering number} of $G$, written $\gamma(G)$, is the minimal cardinality of a normal covering of $G$.
\clearpage

By Corollary~\ref{cor:covering}, if $G$ is a nontrivial finite group, then $\gamma(G) \geq 2$.
When does $\gamma(G)$ achieve this extremal value of $2$? 
There has been much work in this direction, including a complete classification for simple groups by Bubboloni, Spiga and Weigel in 2024 \cite{ref:BubboloniSpigaWeigel24}.
For example, $\gamma(A_n) = 2$ if and only if $4 \leq n \leq 8$, and, in general, $\gamma(A_n)$ grows linearly with $n$, see \cite{ref:EberhardMellon25} for the most recent work in this direction, building on a long line beginning with the work of Bubboloni and Praeger \cite{ref:BubboloniPraeger11}.

We now turn to the second corollary.

\begin{cor}
\label{cor:invariable}
Any subset of a finite group that contains an element from each conjugacy class generates the group.
\end{cor}

This leads to the following definition, introduced by Dixon in 1992 \cite{ref:Dixon92}.
An \emph{invariable generating set} of a group $G$ is a subset $X \subseteq G$ such that for all choices of elements $g_x \in G$ for $x \in X$
\[
G = \langle x^{g_x} \mid x \in X \rangle.
\]

By Corollary~\ref{cor:invariable}, if $G$ is a finite group, then $G$ has an invariable generating set (choose one element from each conjugacy class).
How small an invariable generating set can one find?
For example, famously every finite simple group has a generating set of size two \cite{ref:AschbacherGuralnick84} (for the groups of Lie type, this was proved by Steinberg \cite{ref:Steinberg62}), and one of the many strong versions of this result is that there is, in fact, always an invariable generating set of size two, proved independently in \cite{ref:GuralnickMalle12JLMS} and \cite{ref:KantorLubotzkyShalev11}.

In 2020, Garzoni asked whether some nonabelian finite simple group $G$ has an invariable generating set $\{ g, g^\alpha \}$ for some $g \in G$ and $\alpha \in \Aut(G)$ \cite[Question~2.14]{ref:Garzoni20}, proving that no alternating groups do \cite[Lemma~2.13]{ref:Garzoni20}.
In 2024, Harper \cite[Theorem~3]{ref:Harper24} proved that $\{ g, g^\alpha \}$ is an invariable generating set of the finite simple group $G$ if and only if $g$ is a derangement in every faithful primitive action of the almost simple group $A = \langle G, \alpha \rangle$, and classified all examples.
The examples are $G = \Sp_4(2^f)$ when $A$ contains a graph-field automorphism and $G = \POm^+_d(q)$ with $d = 2^l \geq 8$ when $A$ contains a graph automorphism or a product of graph and field automorphisms.

\vspace{11pt}

We now return to general, possibly infinite, groups.
The class $\mathcal{X}$ of groups satisfying one (and hence all) of the conditions (a)--(c) in Lemma~\ref{lem:abstract} was studied by Wiegold \cite{ref:Wiegold76, ref:Wiegold77}.
Via (a), Theorem~\ref{thm:existence} implies that $\mathcal{X}$ contains all finite groups.
Via (c), it is clear that $\mathcal{X}$ contains all abelian groups.
Via (b), it follows that $\mathcal{X}$ is closed under extensions (exercise). 
In particular, we see that $\mathcal{X}$ includes all virtually soluble groups.
Decades on, there continues to be interesting work in this direction.
For example, Minasyan \cite{ref:Minasyan21} and Goffer and Lazarovich \cite{ref:GofferLazarovich22} independently proved in 2021/2022 that $\mathcal{X}$ is not closed under taking finite index subgroups.

\subsection{Topology: covering spaces}
\label{ss:connections_toplogy}

\sectionquestion{Does every loop lift to a loop in a covering space?}

\begin{ex}
\label{ex:topology}
{\rm
Consider the circle $S^1 = \{ z \in \Com \mid |z| = 1 \}$ and its double cover $p \colon S^1 \to S^1$ defined as $z \mapsto z^2$.
Then the identity map $\phi \colon S^1 \to S^1$ does not lift to the covering space (i.e. there is no continuous map $\widetilde{\phi} \colon S^1 \to S^1$ with $\widetilde{\phi} p = \phi$).
}
\end{ex}

Example~\ref{ex:topology} is typical.
 
\begin{theo}
\label{thm:topology}
Let $X$ be a topological space with a path-connected covering space $p \colon \widetilde{X} \to X$ of finite degree $n > 1$.
Then there exists a continuous map $\phi \colon S^1 \to X$ that does not lift to $\widetilde{X}$.
\end{theo}

\begin{proof}
Fix $x \in X$, let $\Omega = p^{-1}(x)$ and consider the fundamental group $G = \pi_1(X, x)$.
Then $G$ acts on $\Omega$ as follows: for $\alpha \in \Omega$ and $g \in G$, represent $g$ by a loop $g \colon [0, 1] \to X$ with $g(0) = g(1) = x$, lift $g$ to $\widetilde{g} \colon [0, 1] \to \widetilde{X}$ with $\widetilde{g}(0) = \alpha$ and define $\alpha g = \widetilde{g}(1)$, which does not depend on the choice of loop by homotopy lifting.
Since $\widetilde{X}$ is path-connected, this action is transitive.
By Theorem~\ref{thm:existence}, applied to the finite transitive group $G^\Omega$, there is a derangement $g \in G$.
The loop $g$ can be regarded as a continuous map $g \colon S^1 \to X$.
If $g$ lifts to a map $\widetilde{g} \colon S^1 \to \widetilde{X}$, then $\widetilde{g}(0) = \widetilde{g}(1)$, so, by uniqueness of path lifting, $\alpha g = \alpha$ for $\alpha = \widetilde{g}(0)$, which contradicts $g$ being a derangement.
Therefore, $g \colon S^1 \to X$ does not lift to $\widetilde{X}$.
\end{proof}

Our account follows \cite[Section~3]{ref:Serre03}.
\clearpage

\subsection{Computational group theory: finding derangements}
\label{ss:connections_computation}

\sectionquestion{Given generators for a transitive permutation group $G \leq \Sym(n)$, how quickly can we find a derangement?}

Corollary~\ref{cor:proportion} gives a polynomial-time algorithm to find a derangement in $G$: randomly choose elements until you find one!
More precisely, we form a base and strong generating set for $G$ in time polynomial in $n$, with which we produce uniform random elements of $G$ (see \cite[Section~1.13]{ref:Cameron99}, \cite[Section~4.4]{ref:HoltEickOBrien05} or \cite[Section~4]{ref:Lee25}), and we draw elements until one is a derangement.
Each draw succeeds with probability $\delta(G) \geq 1/n$ by Corollary~\ref{cor:proportion}, so the expected number of draws is at most $n$, and more than $cn$ draws are needed with probability at most $(1 - 1/n)^{cn} \leq e^{-c}$. 
Thus, this is a polynomial-time Las Vegas algorithm (see \cite[Section~3.2.1]{ref:HoltEickOBrien05}).

In contrast, in 2017, Arvind \cite{ref:Arvind17} gave a deterministic polynomial-time algorithm, answering a question of Cameron (see \cite[Section~8.1]{ref:CameronWu10} and \cite[Section~4]{ref:Cameron13}).
We follow Cameron's exposition \cite[Section~5]{ref:Cameron13}.
By adapting the proof of the Orbit-Counting Lemma, Arvind proved that for $H \leq G$ and $g \in G$, one can compute
\[
F(Hg) = \frac{1}{|H|} \sum_{h \in H} \fix(hg)
\]
in polynomial time \cite[Lemma~2.4]{ref:Arvind17}.
With this, the algorithm is:
\begin{center}
\begin{tabular}{@{}l@{}}
choose $\omega \in \Omega$ and $g \in G$ with $F(G_\omega g) < 1$, and set $H \coloneqq G_\omega$ \\
\textbf{while} $|Hg| > 1$: \\
\quad choose $\omega \in \Omega$ with $H_\omega < H$ \\
\quad choose $g' \in Hg$ with $F(H_\omega g') \leq F(Hg)$ \\
\quad $H \coloneqq H_\omega$ and $g \coloneqq g'$ \\
\textbf{return} $g$
\end{tabular}
\end{center}

The initial choice is possible as, by the Orbit-Counting Lemma, $F(G) = 1$ is the average of $F(G_\omega g)$ over the cosets $G_\omega g$ while $F(G_\omega) > 1$.
In each subsequent step, $F(Hg)$ is the average of $F(H_\omega g')$ over the cosets $H_\omega g' \subseteq Hg$, so a suitable $g'$ exists.
Since $F(Hg) < 1$ initially and $F(Hg)$ never increases, the element $g$ returned satisfies $\fix(g) = F(\{g\}) < F(G) = 1$, so $g$ is a derangement.

\subsection{Number theory: roots of polynomials modulo primes}
\label{ss:connections_galois}

\sectionquestion{Can a polynomial with no rational root have a root modulo every prime?}

To address this question, write $\Primes$ for the set of primes, and say that a subset $\mathcal{A} \subseteq \Primes$ has \emph{(natural) density} $c$ if
\[
\lim_{N \to \infty} \frac{|\{ x \in \mathcal{A} \mid x \leq N \}|}{|\{ x \in \Primes \mid x \leq N \}|} = c.
\]

\begin{ex}
\label{ex:galois}
{\rm
The polynomial $X^2 + 1$ has no roots over $\Rat$, and it has no root modulo $p \in \Primes$ if and only if $p \equiv 3 \pmod{4}$, a set of primes well known to have density $1/2$.
}
\end{ex}

\begin{theo}
\label{thm:galois}
Let $f \in \Int[X]$ be a monic irreducible polynomial of degree $n > 1$. 
Then there exists a prime $p$ such that $f$ has no roots modulo $p$.
Moreover, the density of such primes is $\delta(\Gal(f)) \geq 1/n$.
\end{theo}

\begin{proof}
Let $\Omega = \{ \omega_1, \dots, \omega_n \}$ be the set of roots of $f$ over $\QQ$, and consider the splitting field $E = \Rat(\omega_1, \dots, \omega_n)$ and Galois group $G = \Gal(E/\Rat)$.
Since $f$ is irreducible, $G$ acts transitively on $\Omega$.
Let $\Primes'$ be the set of primes not dividing the discriminant of $f$, and let $p \in \Primes'$, so $f$ has a set $\Omega_p$ of $n$ distinct roots in $\FF_p$.
The map $\pi_p \colon \lambda \mapsto \lambda^p$ permutes $\Omega_p$ and its orbits correspond to the irreducible factors of $f$ over $\F_p$, so $f$ has no roots in $\F_p$ if and only if $\pi_p$ is a derangement on $\Omega_p$.
To associate to $\pi_p$ an element of $G$, fix a ring homomorphism $\phi \colon \Int[\omega_1, \dots, \omega_n] \to \FF_p$ (via a maximal ideal containing $p$), which restricts to a bijection $\psi \colon \Omega \to \Omega_p$, and define the \emph{Frobenius element} $\Frob_p = \psi \pi_p \psi^{-1} \in \Sym(\Omega)$, which has the same cycle type on $\Omega$ as $\pi_p$ has on $\Omega_p$.
In fact, $\Frob_p$ lies in $G$, and a different choice of $\phi$ gives a conjugate of $\Frob_p$ in $G$, so $\Frob_p$ is well-defined up to conjugacy (see \cite[Section~3.1]{ref:Serre03}).
Since the set of derangements in $G$ is a union of conjugacy classes, Chebotarev's Density Theorem (see \cite{ref:Chebotarev26} for the original statement in terms of analytic density, \cite{ref:Artin24, ref:Artin27} for the version we use in terms of natural density and \cite[Chapter~VII, Section~13]{ref:Neukirch99} for a textbook account) implies that the density of primes $p \in \Primes'$ such that $\Frob_p$ is a derangement equals $\delta(G)$.
By Corollary~\ref{cor:proportion}, $\delta(G) \geq 1/n$.
\end{proof}

\begin{rem}
\label{rem:galois}
{\rm
A typical monic integer polynomial $f$ of degree $n$ has Galois group $\Sym(n)$, in which case, by Theorem~\ref{thm:galois}, the density of primes $p$ for which $f$ has no root is $\delta(\Sym(n)) \approx 1/e$ (see Theorem~\ref{thm:de_montmort}).
More precisely, in 1934, Van~der~Waerden \cite{ref:VanDerWaerden34} proved that for each positive integer $H$, the number of polynomials $X^n + a_1 X^{n-1} + \cdots + a_n \in \Int[X]$ with $|a_i| \leq H$ for all $i$ whose Galois group is not $\Sym(n)$ is $o(H^n)$ as $H \to \infty$.
}
\end{rem}

In 1988, Pollard introduced a new method for factoring integers of a special form, leading to the \emph{(special) number field sieve} (see \cite{ref:Pollard93}).
This algorithm is highly efficient at factoring numbers of the form $a^e \pm b$ where $a$ and $b$ are small, and it was used in \cite{ref:LenstraLenstraManassePollard93} to factor into primes the ninth Fermat number as $2^{2^9}+1 = pqr$ where
\begin{align*}
p &= 2424833 \\
q &= 7455602 825647 884208 337395 736200 454918 783366 342657 \\
r &= 741 640062627530 801524 787141 901937 474059 940781 097519\\ &\qquad\quad  023905 821316 144415 759504 705008 092818 711693 940737.
\end{align*}

The restrictions on the form of the numbers that were amenable to this algorithm motivated the development of the \emph{general number field sieve} (see \cite{ref:BuhlerLenstraPomerance93}).
In these number field sieves, the starting point for factoring a positive integer $N$ is to fix a monic irreducible polynomial $f \in \Int[X]$ of degree $n > 1$ and a positive integer $m$ such that $f(m) \equiv 0 \pmod{N}$.
One challenge that arises in the general number field sieve that does not arise in previous algorithms (including the special number field sieve) is to compute square roots in the ring $\Int[\alpha]$ where $\alpha$ is a root of $f$ in $\QQ$.
This has the potential to dominate the run time of the entire algorithm.
However, a more efficient approach can be taken if one can find an odd prime $p$ such that $f$ has no roots modulo $p$.
This led Lenstra to ask how likely it was to find such a prime (see \cite[Proof of Proposition~9.1]{ref:BuhlerLenstraPomerance93}).
The missing ingredient is precisely Corollary~\ref{cor:proportion}.

Cameron and Cohen write the following in their paper, which appeared in a special issue for Van~Lint's 60th birthday.

\begin{quote}
At the end of May, 1991, H. W. Lenstra Jr asked the second author whether the number of fixed point free elements of a finite transitive group on a set of $n$ elements, $n > 1$, is at least $1/n$th of the order of~$G$ (with equality only for
Frobenius groups of order $n(n - 1)$). 
An affirmative answer would be useful to know in connection with a new number factorization algorithm, so he was told.
Almost convinced that such a result should be in the pre-world-war literature and taking a lazy option, the second author e-mailed a few colleagues for references. 
In the meantime, he tried some of the combinatorial approaches taught by van~Lint. 
The result was the proposition below, which generalizes an affirmative answer to the above query to arbitrary $n \geq 1$. 
[\dots]
Before reading the solution (e-mailed shortly after the reference query), the
first author replied with more or less the same result. No surprise, one might
argue, for someone who has also been under van Lint’s influence.
\end{quote}

For more on this topic, see \cite{ref:Serre03}.

\subsection{Number theory: values of polynomials over finite fields}
\label{ss:connections_values}

\sectionquestion{How many values can a polynomial take on a finite field?}

To address this question, we use a function field analogue of Section~\ref{ss:connections_galois}, in which $\Int$ is replaced by $\F_q[t]$ and reduction modulo a prime by specialising $t$ to an element $a \in \F_q$.
Let $q$ be a power of a prime $p$, let $f \in \F_q[X]$ have degree $n$ and let $V_f$ be the number of values that $f$ takes on $\F_q$.
As $x \mapsto x^p$ is a bijection of $\F_q$, replacing $f(X)$ by $f(X^p)$ does not change $V_f$, so we assume that $f$ is not a polynomial in $X^p$, and, as polynomials of degree~$1$ are bijections, we assume that $n \geq 2$.

\begin{ex}
\label{ex:values}
{\rm
Let $f = X^n$ with $p \nmid n$, and let $d = \gcd(n, q-1)$.
The nonzero values of $f$ form the subgroup of index $d$ in the cyclic group $\F_q^\times$, so $V_f = 1 + (q-1)/d$.
For example, if $q$ is odd, then $X^2$ takes $(q+1)/2$ values, while $X^3$ takes every value in $\F_5$.
}
\end{ex}

An element $a \in \F_q$ is a value of $f$ exactly when $f(X) - a$ has a root in $\F_q$, so, for an indeterminate $t$, $1 - V_f/q$ is the proportion of $a \in \F_q$ for which $f(X) - t$ has no root in $\F_q$ upon specialising $t \mapsto a$, analogous to the density in Theorem~\ref{thm:galois}.

Let $E$ be the splitting field of $f(X) - t$ over $\F_q(t)$, let $\Omega$ be the set of roots of $f(X) - t$ and let $G$ and $N$ be the Galois groups of $f(X) - t$ over $\F_q(t)$ and $\FF_q(t)$, acting on $\Omega$, the \emph{arithmetic} and \emph{geometric monodromy groups} of $f$.
The elements of $E$ that are algebraic over $\F_q$ form a finite field $\F_{q^r}$, and restriction to $\F_{q^r}$ gives a surjection from $G$ onto the cyclic group $\Gal(\F_{q^r}/\F_q)$ with kernel $N$, so the preimage of the Frobenius map $\lambda \mapsto \lambda^q$ is a coset $Nx$, the \emph{Frobenius coset}.
Guralnick and Wan proved the following \cite[Section~1]{ref:GuralnickWan97}.

\begin{theo}
\label{thm:values}
Let $f \in \F_q[X]$ be a polynomial of degree $n \geq 2$ that is not a polynomial in $X^p$.
Let $G$ and $N$ be the arithmetic and geometric monodromy groups of $f$ and let $Nx$ be the Frobenius coset.
Then
\[
V_f = (1 - \delta(Nx))q + O(\sqrt{q}).
\]
Moreover, if $G = N$, then $V_f \leq (1-1/n)q + O(\sqrt{q})$.
(The implied constants depend only on $n$.)
\end{theo}

\begin{proof}[Proof sketch]
For $k = \F_q$ or $k = \FF_q$, the polynomial $f(X) - t$ is irreducible in $k[X][t] = k[t][X]$ and primitive over $k[t]$, so it is irreducible over $k(t)$ by Gauss's Lemma, and hence $G$ and $N$ act transitively on $\Omega$.
Moreover, $f(X) - t$ is separable, since its derivative $f'(X)$ is nonzero as $f$ is not a polynomial in $X^p$.
Hence, $|\Omega| = n$.

Let $a \in \F_q$.
The map $\lambda \mapsto \lambda^q$ permutes the roots of $f(X) - a$ in $\FF_q$, and its fixed points are exactly the roots in $\F_q$, so $a$ is a value of $f$ if and only if this permutation has a fixed point.
As for the Frobenius element in the proof of Theorem~\ref{thm:galois}, for all but boundedly many $a$, specialising $t$ to $a$ gives an element $\Frob_a \in G$, well-defined up to conjugacy, which has the same cycle type on $\Omega$ as this permutation has on the roots of $f(X) - a$, and $\Frob_a \in Nx$ since it acts on $\F_{q^r}$ as $\lambda \mapsto \lambda^q$.
The function field analogue of Chebotarev's Density Theorem gives $V_f = (1 - \delta(Nx))q + O(\sqrt{q})$.
If $G = N$, then, by Corollary~\ref{cor:proportion}, $\delta(Nx) = \delta(N) \geq 1/n$.
\end{proof}

\begin{ex}
\label{ex:values_galois}
{\rm
Let $f = X^n$ as in Example~\ref{ex:values}.
The roots of $X^n - t$ are $\zeta^i\alpha$ for $i \in \Int/n\Int$, where $\alpha^n = t$ and $\zeta \in \FF_q$ is a primitive $n$th root of unity, so we identify $\Omega$ with $\Int/n\Int$.
Then $N$ consists of the translations $i \mapsto i + j$ and, as every element of $Nx$ maps $\zeta$ to $\zeta^q$, the coset $Nx$ consists of the maps $i \mapsto qi + j$, which fix a point if and only if $d$ divides $j$.
Hence, $\delta(Nx) = 1 - 1/d$, and Theorem~\ref{thm:values} gives $V_f = q/d + O(\sqrt{q})$, agreeing with Example~\ref{ex:values}.
}
\end{ex}

\begin{rem}
\label{rem:values}
{\rm
Like in Remark~\ref{rem:galois}, a typical $f \in \F_q[X]$ satisfies $G = N = \Sym(n)$ (see \cite{ref:Cohen70}), so ${V_f = (1 - \delta(\Sym(n)))q + O(\sqrt{q})}$, which is roughly $(1 - 1/e)q$, by Theorem~\ref{thm:values}.
This was proved by Birch and Swinnerton-Dyer in 1959 \cite{ref:BirchSwinnertonDyer59}, answering a question of Chowla \cite{ref:Chowla52}.
}
\end{rem}

In 1993, Wan \cite{ref:Wan93} proved that if $V_f < q$, then $V_f \leq q - (q-1)/n$, which is sharp (see \cite[Section~1]{ref:GuralnickWan97}).
This motivated Guralnick and Wan \cite[Theorem~1.1]{ref:GuralnickWan97} to find the next possible value of $\delta(Nx)$ after $1/n$. 
The case $G = N$ is Theorem~\ref{thm:guralnick_wan}, which gives ${V_f \leq (1 - 2/n)q + O(\sqrt{q})}$ for $n \geq 10$ unless $N$ is a Frobenius group of order $n(n-1)$.
In general, they \cite[Theorem~1.2]{ref:GuralnickWan97} proved that if $p \nmid n$ and $\delta(Nx) > 0$, then $\delta(Nx) > 1/6$, so ${V_f \leq 5q/6 + O(\sqrt{q})}$.
Theorem~\ref{thm:garzoni} gives ${V_f \leq (1 - (\sqrt{n}+1)/(2n))q + O(\sqrt{q})}$ for $G = N$ and sufficiently large $n$, unless $N$ is a primitive Frobenius group.

\subsection{Number theory: Kronecker equivalence}
\label{ss:connections_kronecker}

\sectionquestion{If two irreducible polynomials have roots modulo exactly the same primes, must they have the same degree?}

Fix a number field $k$, i.e. a finite extension of $\Rat$.
The \emph{primes} of $k$ are the nonzero prime ideals of its ring of integers $\mathcal{O}_k$, which for $k = \Rat$ correspond to the usual primes.
For a finite extension $K/k$, the \emph{Kronecker set} $D(K/k)$ is the set of primes $\mathfrak{p}$ of $k$ for which some prime $\mathfrak{P}$ of $K$ containing $\mathfrak{p}$ satisfies $\mathcal{O}_K/\mathfrak{P} = \mathcal{O}_k/\mathfrak{p}$.
For $k = \Rat$ and $K = \Rat(\alpha)$, a prime $p$ not dividing the discriminant of the minimal polynomial $f$ of $\alpha$ lies in $D(K/\Rat)$ if and only if $f$ has a root modulo $p$.
Two extensions $K$ and $L$ of $k$ are \emph{Kronecker equivalent} if $D(K/k)$ and $D(L/k)$ differ by only finitely many primes.

\begin{ex}
\label{ex:kronecker}
{\rm
Let $f = X^6 - 3X^2 - 1$ and $g = X^3 - 3X - 1$, which are irreducible, and let $\alpha$ be a root of $f$, so $\alpha^2$ is a root of $g$.
Let $p \neq 3$ be prime.
Clearly, if $f$ has a root modulo $p$, then so does $g$, and we prove the converse (showing that each has a root modulo $p$ if and only if $p \equiv \pm 1 \pmod{9}$).
Let $\zeta \in \FF_p$ be a primitive 9th root of unity, so $-(\zeta^{k} + \zeta^{-k})$ is a root of $g$ for $k \in \{1, 2, 4\}$.
The map $\lambda \mapsto \lambda^p$ sends $-(\zeta^k + \zeta^{-k})$ to $-(\zeta^{pk} + \zeta^{-pk})$, so one root of $g$ is in $\F_p$ if and only if they all are (if and only if $p \equiv \pm 1 \pmod{9}$).
Hence, if $g$ has a root modulo $p$, then it has three, whose product is $1$, so one of them is a square in $\F_p$, and its square roots are roots of $f$.
Therefore, $f$ has a root modulo $p$ if and only if $g$ does, so $K = \Rat(\alpha)$ and $L = \Rat(\alpha^2)$ are Kronecker equivalent.
}
\end{ex}

In 1977, Jehne \cite{ref:Jehne77} connected this with derangements.

\begin{lem}
\label{lem:jehne}
Let $K$ and $L$ be finite extensions of $k$ contained in a Galois extension $E$ of $k$, and let $A = \Gal(E/k)$, $U = \Gal(E/K)$ and $V = \Gal(E/L)$.
Then the following are equivalent
\begin{itemize}
\item[{\rm (a)}] $K$ and $L$ are Kronecker equivalent
\item[{\rm (b)}] $\bigcup_{a \in A} U^a = \bigcup_{a \in A} V^a$
\item[{\rm (c)}] the same elements of $A$ are derangements on $A/U$ and $A/V$.
\end{itemize}
\end{lem}

\begin{proof}[Proof sketch when $k = \Rat$]
Clearly (b) $\Leftrightarrow$ (c).
For (a) $\Leftrightarrow$ (c), write $K = \Rat(\alpha)$ and $L = \Rat(\beta)$, and let $f$ and $g$ be the minimal polynomials of $\alpha$ and $\beta$, so the actions of $A$ on their roots are equivalent to its actions on $A/U$ and $A/V$.
We apply the construction of the Frobenius element in the proof of Theorem~\ref{thm:galois} to the minimal polynomial $h$ of $\theta$, where $E = \Rat(\theta)$.
For all but finitely many primes $p$, this gives $\Frob_p \in A$, well-defined up to conjugacy, acting on the roots of $h$, and hence, as $\alpha$ and $\beta$ are rational polynomials in $\theta$, on those of $f$ and $g$, as $\lambda \mapsto \lambda^p$ acts on their reductions modulo $p$.
Hence, $p \in D(K/\Rat)$ if and only if $\Frob_p$ is not a derangement on $A/U$, and similarly for $L$ and $V$, which proves (c) $\Rightarrow$ (a).
For \textsc{not} (c) $\Rightarrow$ \textsc{not} (a), if, say, $x \in A$ is a derangement on $A/V$ but not on $A/U$, then, by Chebotarev's Density Theorem (see the proof of Theorem~\ref{thm:galois}), $x$ is conjugate to $\Frob_p$ for infinitely many $p$, giving infinitely many primes in $D(K/\Rat) \setminus D(L/\Rat)$.
\end{proof}

\begin{ex}
\label{ex:jehne}
{\rm
We return to Example~\ref{ex:kronecker}.
Let $E$ be the splitting field of $f$, and let $A = \Gal(E/\Rat)$ acting on the set of roots of $f$.
One can verify that $A \cong A_4$, acting on the six $2$-subsets of $\{1, 2, 3, 4\}$ as in Example~\ref{ex:k-subsets}.
Then $U = \Gal(E/K)$ corresponds to the stabiliser of a $2$-subset, say $U = \langle (1 \ 2)(3 \ 4) \rangle$, and $V = \Gal(E/L)$ contains $U$ with index $[K:L] = 2$, so $V = \langle (1 \ 2)(3 \ 4), (1 \ 3)(2 \ 4) \rangle$.
Here $\bigcup_{a \in A} U^a = V = \bigcup_{a \in A} V^a$, so Lemma~\ref{lem:jehne} implies that $K$ and $L$ are Kronecker equivalent, as shown directly in Example~\ref{ex:kronecker}.
Moreover, for all but finitely many primes $p$, the polynomials $f$ and $g$ have a root modulo $p$ exactly when the Frobenius element $\Frob_p$ fixes a root of $f$, equivalently of $g$, that is, when $\Frob_p \in V$, so the density of such primes is $|V|/|A| = 1/3$, as the condition $p \equiv \pm 1 \pmod 9$ in Example~\ref{ex:kronecker} predicts.
}
\end{ex}

If $K$ is a quadratic extension of $k$ and $L$ is a finite extension of $k$ that is Kronecker equivalent to $K$, then, in the notation of Lemma~\ref{lem:jehne}, $U \trianglelefteqslant A$ has index $2$ and $\{V, V^a\}$ is a normal covering of $U$ for $a \in A \setminus U$.
Is this possible when $V < U$?
Jehne reduced this question to the case where $U$ is simple, and in 1988, Saxl \cite{ref:Saxl88} handled this case using the CFSG, proving the following.

\begin{theo}
\label{thm:saxl}
If $K/k$ is a quadratic extension of number fields, then every extension of $k$ that is Kronecker equivalent to $K$ is isomorphic to $K$.
\end{theo}

More generally, Neumann and Praeger conjectured the following (see \cite[Conjectures~4.1 and~4.1$'$]{ref:Praeger94}).

\begin{conj}
\label{conj:neumann_praeger}
There exists a function $h \colon \Nat \to \Nat$ such that the following hold.
\begin{itemize}
\item [{\rm (a)}]
If $K/k$ and $L/k$ are Kronecker equivalent extensions of number fields with $[K:k] = n$, then $[L:k] \leq h(n)$.
\item [{\rm (b)}]
If $A$ is a finite group with subgroups $U$ and $V$ with $|A:U| = n$ and $\bigcup_{a \in A} U^a = \bigcup_{a \in A} V^a$, then $|A:V| \leq h(n)$.
\item [{\rm (c)}]
If $A$ is a finite group, $V \leq U \trianglelefteqslant A$ and $\{V^{a_1}, \dots, V^{a_n}\}$ is a normal covering of $U$ for coset representatives $a_1, \dots, a_n$ of $U$ in $A$, then $|U:V| \leq h(n)$.
\end{itemize}
\end{conj}

The three parts are equivalent, up to changing the function $h$.
Indeed, (a) and (b) are equivalent by Lemma~\ref{lem:jehne} and Remark~\ref{rem:brauer} below, (b) implies (c) since $\{V^{a_1}, \dots, V^{a_n}\}$ is a normal covering of $U$ exactly when $U = \bigcup_{a \in A} V^a$, and Praeger showed that (c) implies (b) by passing to the core of $U$ \cite[p.~26]{ref:Praeger94}.

\subsection{Number theory: central simple algebras}
\label{ss:connections_brauer}

\sectionquestion{How many central simple algebras over $\Rat$ become matrix algebras over a given finite extension?}

Fix a field $k$.
A \emph{central simple algebra} over~$k$ is a finite-dimensional associative algebra over~$k$ that is simple as a ring and whose centre is $k$.
By Wedderburn's Theorem \cite{ref:Wedderburn08}, every central simple algebra over~$k$ is isomorphic to a matrix algebra $M_d(D)$ where $d$ is a positive integer and $D$ is a division algebra over~$k$.
We say that central simple algebras $A$ and $B$ over~$k$ are \emph{equivalent} if $A \cong M_d(D)$ and $B \cong M_e(D)$ for the same division algebra $D$.
More generally, for an extension $K/k$, central simple algebras $A$ and $B$ over~$k$ are \emph{equivalent over~$K$} if the tensor products $A \otimes_k K$ and $B \otimes_k K$ are equivalent as central simple algebras over~$K$.

\begin{ex}
\label{ex:brauer}
{\rm
Two examples of central simple algebras over~$\Real$ are the real numbers $\Real$ itself and the quaternions $\Ham$.
These are, in fact, both division rings over~$\Real$, so the form guaranteed by Wedderburn's Theorem in this case is simply $\Real \cong M_1(\Real)$ and $\Ham \cong M_1(\Ham)$.
The other finite-dimensional division algebra over $\Real$ is the complex numbers $\Com$, but this is not a central simple algebra over~$\Real$ since its centre is $\Com$.
In view of this, $\Real$ and $\Ham$ represent the two distinct equivalence classes of central simple algebras over $\Real$.

Now consider the corresponding algebras over~$\Com$ via the tensor product $A \mapsto A \otimes_{\Real} \Com$.
Certainly $\Real \otimes_{\Real} \Com = \Com \cong M_1(\Com)$ and we see that $\Ham \otimes_{\Real} \Com = M_2(\Com)$ by considering the representation
\[
\Ham \to M_2(\Com) \quad a + bi + cj + dk \mapsto \left( \begin{array}{cc} a + bi & c + di \\ -c + di & a - bi \end{array} \right).
\]
Therefore, all central simple algebras over~$\Real$ are equivalent over~$\Com$.
}
\end{ex}

In contrast, the behaviour is very different for finite extensions of \emph{global fields}, that is, finite extensions of $\Rat$ (\emph{number fields}) or of $\F_q(t)$ (\emph{function fields} in one variable over a finite field).

\begin{theo}
\label{thm:brauer}
Let $K/k$ be a finite nontrivial extension of global fields.
Then there are infinitely many inequivalent central simple algebras over $k$ that are equivalent to $k$ over $K$.
\end{theo}

To prove this, write $[A]$ for the equivalence class of a central simple algebra $A$ over $k$, and $\mathrm{Br}(k)$ for the set of such equivalence classes.
One can show that $[A] \cdot [B] = [A \otimes_k B]$ is a well-defined commutative binary operation on $\mathrm{Br}(k)$, with $[k]$ as the identity and $[A^{\rm op}]$ as the inverse of $[A]$.
This makes $\mathrm{Br}(k)$ an abelian group, known as the \emph{Brauer group} of~$k$.
Moreover, the map $[A] \mapsto [A \otimes_k K]$ is a homomorphism $\mathrm{Br}(k) \to \mathrm{Br}(K)$ whose kernel is the \emph{relative Brauer group} $\mathrm{Br}(K/k)$, consisting of the classes of central simple algebras over $k$ that are equivalent to $k$ over $K$.

When $k = \Rat$, the Brauer group is described by the classification of central simple algebras over $\Rat$, due to Albert, Brauer, Hasse and Noether (see \cite[Section~32]{ref:Reiner75}).
To each class $[A]$ one attaches an \emph{invariant} $\mathrm{inv}_p[A] \in \Rat/\Int$ for each prime $p$, together with an invariant $\mathrm{inv}_\infty[A] \in \{0, 1/2 \}$ recording whether $A$ is equivalent to $\Real$ or to $\Ham$ over $\Real$.
The invariants are additive, only finitely many of them are nonzero, they sum to $0$, they determine $[A]$ and every such family arises.
In other words, there is an exact sequence
\[
0 \to \mathrm{Br}(\Rat) \xrightarrow{\rho} \tfrac{1}{2}\Int/\Int \oplus \bigoplus_{p} \Rat/\Int \xrightarrow{\sigma} \Rat/\Int \to 0,
\]
where $\rho$ records the invariants and $\sigma$ sums them.

\begin{proof}[Proof sketch of Theorem~\ref{thm:brauer} when $k = \Rat$]
We show that $\mathrm{Br}(K/\Rat)$ is infinite.
Write $K = \Rat(\alpha)$ for an algebraic integer $\alpha$, let $f \in \Int[X]$ be the minimal polynomial of $\alpha$, let $\Omega$ be the set of roots of $f$ in its splitting field $E$ and let $G = \Gal(E/\Rat)$, which acts transitively on $\Omega$.
For a prime $p$ not dividing the discriminant $\mathrm{disc}(f)$, let $n_1, \dots, n_l$ be the degrees of the irreducible factors of $f$ modulo $p$, which are the cycle lengths on $\Omega$ of ${\Frob_p \in G}$ (see the proof of Theorem~\ref{thm:galois}), and let $m_p = \gcd(n_1, \dots, n_l)$.
Moreover, $[A] \in \mathrm{Br}(K/\Rat)$ if and only if $m_p \, \mathrm{inv}_p[A] = 0$ for every prime $p \nmid \mathrm{disc}(f)$ and analogous conditions hold at the primes $p \mid \mathrm{disc}(f)$ and at $p = \infty$ (see the proof of \cite[Theorem~2]{ref:FeinKantorSchacher81}).

By Theorem~\ref{thm:fein_kantor_schacher}, $G$ has a derangement $x$ of $r$-power order for some prime $r$, so every cycle of $x$ on $\Omega$ has length divisible by $r$.
By Chebotarev's Density Theorem (see the proof of Theorem~\ref{thm:galois}), $x$ is conjugate to $\Frob_p$ for infinitely many primes $p_0, p_1, p_2, \dots$ not dividing $\mathrm{disc}(f)$, so $r$ divides $m_{p_i}$ for each $i$.
Therefore, for each $i \geq 1$, the family with $1/r$ at $p_0$, $-1/r$ at $p_i$ and $0$ elsewhere lies in $\ker \sigma = \im \rho$, so it is $\rho([A_i])$ for some $[A_i] \in \mathrm{Br}(\Rat)$, which lies in $\mathrm{Br}(K/\Rat)$.
These classes are distinct, so $\mathrm{Br}(K/\Rat)$ is infinite.
\end{proof}

The above proof can be reversed.
Indeed, assume that $G$ has no derangement of prime-power order.
Then $m_p = 1$ for every prime $p \nmid \mathrm{disc}(f)$, as otherwise, for a prime $r$ dividing $m_p$, raising $\Frob_p$ to the $r'$-part of its order would give a derangement of $r$-power order (see Remark~\ref{rem:p-elements}).
Hence, the invariants of each element of $\mathrm{Br}(K/\Rat)$ vanish at all but the finitely many primes dividing $\mathrm{disc}(f)$ and they lie in the finite group $\frac{1}{n}\Int/\Int$ as each such element has order dividing $n = [K : \Rat]$ \cite[Theorems~28.5 and~29.22]{ref:Reiner75}, so there are finitely many choices of invariants and $\mathrm{Br}(K/\Rat)$ is finite.

\begin{ex}
\label{ex:brauer_gaussian}
{\rm
Let $K = \Rat(i)$, the splitting field of $f = X^2 + 1$.
For each odd prime $p$, by Example~\ref{ex:galois}, $f$ is irreducible modulo $p$ (equivalently, $m_p > 1$) if and only if $p \equiv 3 \pmod 4$.
For nonzero $a, b \in \Rat$, the \emph{quaternion algebra} $(a,b)_\Rat$ has basis $1, i, j, ij$ with $i^2 = a$, $j^2 = b$ and $ij = -ji$, and its invariant at a prime $\ell$ or at $\infty$ is $0$ if $ax^2 + by^2 = z^2$ has a nonzero solution over $\Rat_\ell$ or $\Real$, respectively, and $1/2$ otherwise.
Fix a prime $p > 3$ with $p \equiv 3 \pmod 4$, and let $A_p = (-1, 3p)_\Rat$.
As $-1$ is not a square modulo $3$ or $p$, the conic $-x^2 + 3py^2 = z^2$ has no nonzero solution over $\Rat_3$ or $\Rat_p$, but it has one over $\Real$ and over $\Rat_\ell$ for each odd prime $\ell \notin \{3, p\}$, so, as the invariants sum to $0$, the invariants of $A_p$ are $1/2$ at $3$ and $p$ and $0$ elsewhere.
Moreover, $A_p$ is $4$-dimensional over $\Rat$ and contains $\Rat(i)$, so $A_p \otimes_\Rat \Rat(i) \cong M_2(\Rat(i))$.
Hence, the algebras $A_p$ are pairwise inequivalent over $\Rat$ but all equivalent over $\Rat(i)$.
}
\end{ex}

\begin{rem}
\label{rem:brauer}
{\rm
For $k = \Rat$, the argument above shows that $\mathrm{Br}(K/k)$ is infinite if and only if the Galois group has a derangement of prime-power order, and Fein, Kantor and Schacher \cite{ref:FeinKantorSchacher81} showed that the same holds for every global field $k$.
Moreover, every finite transitive permutation group $G \leq \Sym(\Omega)$ with $|\Omega| > 1$ arises as such a Galois group over some global field $k$, so Theorem~\ref{thm:brauer} is equivalent (via a CFSG-free proof) to Theorem~\ref{thm:fein_kantor_schacher}, as observed in \cite[p.~48]{ref:FeinKantorSchacher81}.
}
\end{rem}

The systematic study of this topic began with Fein and Schacher's paper \cite{ref:FeinSchacher81}.
Later in 1981, Fein, Kantor and Schacher \cite{ref:FeinKantorSchacher81} proved Theorem~\ref{thm:brauer} by establishing the required result on derangements, namely Theorem~\ref{thm:fein_kantor_schacher}, whose proof uses Kantor's classification of the linear groups generated by transvections \cite{ref:Kantor79}.
Our account benefited from \cite[Section~4.1]{ref:Cameron11}.

\subsection{Representation theory: vanishing of characters}
\label{ss:connections_character_vanishing}

\sectionquestion{Which characters vanish on some element?}

For this section, fix a nontrivial finite group $G$ and a proper subgroup $H$ of $G$.
The induced character $1_H^G$ is the permutation character of the action of $G$ on $G/H$, so $1_H^G(g) = \fix(g)$ for all $g \in G$, and, by Theorem~\ref{thm:existence}, $1_H^G$ \emph{vanishes} on some $g \in G$, that is, $1_H^G(g) = 0$.
More generally, if $\chi$ is a character of $H$, then
\[
\chi^G(g) = \frac{1}{|H|} \sum_{x \in G,\, g^x \in H} \chi(g^x),
\]
so $\chi^G$ vanishes on all derangements in $G$ on $G/H$.
Hence, every character induced from a proper subgroup vanishes on some element.

For irreducible characters, in 1904, Burnside \cite{ref:Burnside04} proved the following (see also \cite[Theorem~3.15]{ref:Isaacs76}).

\begin{theo}
\label{thm:characters_burnside}
Let $G$ be a finite group with a nonlinear irreducible character $\chi$.
Then $\chi(g) = 0$ for some $g \in G$.
\end{theo}

We first record a preliminary result (see \cite[Lemma~3.14]{ref:Isaacs76}).

\begin{lem}
\label{lem:characters_cyclic}
Let $C$ be a finite cyclic group, let $S = \{ x \in C \mid \langle x \rangle = C \}$ and let $\chi$ be a character of $C$.
Then
\begin{itemize}
\item[{\rm (a)}]  $\prod_{x \in S} \chi(x)$ and $\sum_{x \in S} \chi(x)$ are integers
\item[{\rm (b)}] if $\chi(x) \neq 0$ for all $x \in S$, then 
\[
\prod_{x \in S} |\chi(x)|^2 \geq 1 \quad \text{and} \quad \sum_{x \in S} |\chi(x)|^2 \geq |S|.
\]
\end{itemize}
\end{lem}

\begin{proof}
Write $n = |C|$, let $E$ be the splitting field of $X^n-1$ over $\Rat$ and let $\sigma \in \Gal(E/\Rat)$.
Then there exists an integer $m$ coprime to $n$ such that if $\epsilon^n = 1$, then $\epsilon^\sigma = \epsilon^m$.
For all $x \in C$, we can write $\chi(x) = \epsilon_1 + \cdots + \epsilon_k$ where $\epsilon_1^n = \cdots = \epsilon_k^n = 1$, so 
\[
\chi(x)^\sigma = \epsilon_1^\sigma + \cdots + \epsilon_k^\sigma = \epsilon_1^m + \cdots + \epsilon_k^m = \chi(x^m).
\]
Since $m$ is coprime to $n$, the set $S$ is invariant under $x \mapsto x^m$, so $\prod_{x \in S} \chi(x)$ and $\sum_{x \in S} \chi(x)$ are invariant under $\sigma$, which forces these algebraic integers to be (rational) integers, as claimed in (a).

Now assume that $\chi(x) \neq 0$ for all $x \in S$.
Since $|\chi(x)|^2 = (\chi\overline{\chi})(x)$ and $\chi\overline{\chi}$ is a character of $C$, part~(a) implies that $\prod_{x \in S} |\chi(x)|^2$ is an integer.
Moreover, $\prod_{x \in S} |\chi(x)|^2 \geq 1$ as $|\chi(x)|^2 > 0$ for all $x \in S$.
By the inequality of arithmetic and geometric means,
\[
\frac{1}{|S|} \sum_{x \in S} |\chi(x)|^2 \geq \Bigl( \prod_{x \in S} |\chi(x)|^2 \Bigr)^{1/|S|} \geq 1,
\]
completing the proof of (b).
\end{proof}

We can now prove Theorem~\ref{thm:characters_burnside}.

\begin{proof}[Proof of Theorem~\ref{thm:characters_burnside}]
For a contradiction, suppose that $\chi(g) \neq 0$ for all $g \in G$.
Partition $G$ into equivalence classes, where $x, y \in G$ are equivalent if $\langle x \rangle = \langle y \rangle$.
For each class $S$, Lemma~\ref{lem:characters_cyclic}(b) gives $\sum_{x \in S}|\chi(x)|^2 \geq |S|$.
Hence, $\sum_{x \in G \setminus \{1\}}|\chi(x)|^2 \geq |G| - 1$, so $|G| = \sum_{x \in G} |\chi(x)|^2 \geq |G| - 1 + \chi(1)^2$, which implies that $\chi(1) = 1$, which contradicts our assumption that $\chi$ is nonlinear.  
\end{proof}

For a character $\chi$ of a proper subgroup $H < G$, we know that $\chi^G$ vanishes on any derangement in $G$ on $G/H$, so Theorem~\ref{thm:fein_kantor_schacher} actually guarantees that $\chi^G$ vanishes on an element of prime-power order.
This suggests a stronger version of Theorem~\ref{thm:characters_burnside}, which, in 2000, Malle, Navarro and Olsson \cite{ref:MalleNavarroOlsson00} proved.

\begin{theo}
\label{thm:malle_navarro_olsson}
Let $G$ be a finite group with a nonlinear irreducible character $\chi$.
Then $\chi(g) = 0$ for some $g \in G$ of prime-power order.
\end{theo}

If $\chi$ is induced from a proper subgroup of $G$, then $\chi$ vanishes on a prime-power order element as an immediate consequence of Theorem~\ref{thm:fein_kantor_schacher}, so one reduces to the case where $\chi$ is not induced from a proper subgroup, and this is the first step involved in a more general reduction to the case where $G$ is simple.
With this reduction, the proof then invokes the CFSG explicitly before employing tools from the classical representation theory of the symmetric group and Lusztig's classification of the irreducible complex characters of the groups of Lie type, together with the known character tables of the sporadic groups.
Note that in this proof the reduction to simple groups itself uses the CFSG, indirectly via Theorem~\ref{thm:fein_kantor_schacher}.

\subsection{Representation theory: restriction of characters}
\label{ss:connections_character_restriction}

\sectionquestion{Can a proper subgroup distinguish all characters of a group?}

As in the previous section, fix a nontrivial finite group $G$.
By Corollary~\ref{cor:covering}, not every element of $G$ is conjugate into a given proper subgroup.
In 2017, Isaacs and Navarro \cite{ref:IsaacsNavarro17} used this to give an analogue for characters.
The result is elementary, and Serre had asserted it in 2003 \cite[Theorem~4'']{ref:Serre03}.

\begin{theo}
\label{thm:characters_restriction}
Let $G$ be a finite group and let $H$ be a proper subgroup of $G$.
Then there exist distinct characters $\psi$ and $\phi$ of $G$ such that the restrictions $\psi_H$ and $\phi_H$ are equal.
\end{theo}

We first prove a lemma (see \cite[Lemma~B]{ref:IsaacsNavarro17}).

\begin{lem}
\label{lem:characters_sum}
Let $C$ be a finite cyclic group and let $S = \{ x \in C \mid \langle x \rangle = C \}$.
Then there exists a character $\chi$ of $C$ such that $|C|$ divides $\chi(1)$ and for all nontrivial elements $x \in C$
\begin{equation}
\label{eq:chi}
\chi(x) = \left\{ \begin{array}{ll} |C| & \text{if $x \in S$} \\ 0 & \text{if $x \not\in S$.} \end{array} \right.
\end{equation}
\end{lem}

\begin{proof}
Define $\theta\colon C \to \Com$ as 
\[
\theta(x) = \left\{ \begin{array}{ll} |C| & \text{if $x \in S$} \\ 0 & \text{if $x \not\in S$.} \end{array} \right.
\]
Since $C$ is abelian, $\theta$ is a class function, so $\theta = \sum_{\lambda \in \text{Irr}(C)} a_\lambda \lambda$ where
\[
a_\lambda = \frac{1}{|C|}\sum_{x \in C}\theta(x)\overline{\lambda(x)} = \sum_{x \in S}\overline{\lambda(x)},
\]
which is an integer by Lemma~\ref{lem:characters_cyclic}(a).
Hence, $\theta$ is a generalised character, so for the regular character $\rho$ of $C$ and a large enough integer $m$, the sum $\chi = \theta + m \rho$ is a character.
Moreover, $\chi(x) = \theta(x)$ for all nontrivial $x \in C$, and $\chi(1) = \theta(1) + m|C|$ is divisible by $|C|$, since $\theta$ only takes the values $0$ and $|C|$.
\end{proof}

We can now prove Theorem~\ref{thm:characters_restriction}.

\begin{proof}[Proof of Theorem~\ref{thm:characters_restriction}]
By Corollary~\ref{cor:covering}, we can fix an element $c \in G$ that is not conjugate into $H$.
Write $C = \langle c \rangle$ and $S = \{ x \in C \mid \langle x \rangle = C \}$, and, by Lemma~\ref{lem:characters_sum}, fix a character $\chi$ of $C$ satisfying \eqref{eq:chi} for all nontrivial $x \in C$.
Let $h$ be a nontrivial element of $H$.
If a conjugate of $h$ is contained in $S$, then a generator of $C$ (and hence $c$) is conjugate to an element of $H$, which contradicts the construction.
Therefore, no conjugate of $h$ is contained in $S$, so $\chi^G$ vanishes on all nonidentity elements of $H$.
The regular character $\rho$ of $G$ also vanishes on all nonidentity elements. 
Let $\phi = \chi^G$ and $\psi = (\chi(1)/|C|)\rho$, which are characters of $G$ with $\phi(1) = |G:C|\chi(1) = \psi(1)$.
Therefore, $\phi_H = \psi_H$.
However, $\chi^G(c)$ is nonnegative and the term for $x = 1$ is $\chi(c) = |C|$, so $\phi(c) \geq 1 > 0 = \psi(c)$ and $\phi \neq \psi$.
\end{proof}

\subsection{Combinatorics: intersecting families}
\label{ss:connections_intersecting}

\sectionquestion{How large can a family of permutations be if any two of them agree on some point?}

For this section, fix a finite set $\Omega$ of size $n$.
A subset $S$ of the power set $\P(\Omega)$ is an \emph{intersecting family of subsets} if any two subsets in $S$ intersect  (i.e. for all $X, Y \in S$ there exists $\omega \in \Omega$ such that $\omega \in X \cap Y$).
The Erd\H{o}s--Ko--Rado Theorem \cite{ref:ErdosKoRado61} states that if $S$ is an intersecting family of $k$-subsets of $\Omega$, where $n > 2k$, then $|S| \leq \tbinom{n-1}{k-1}$, with equality if and only if $S$ consists of all $k$-subsets of $\Omega$ containing a given point.

This classical theorem has motivated much work in extremal combinatorics (see, for example, the survey \cite{ref:Ellis22}).
By analogy, a subset $S$ of the symmetric group $\Sym(\Omega)$ is an \emph{intersecting family of permutations} if any two permutations in $S$ agree on at least one point (i.e. for all $x, y \in S$ there exists $\omega \in \Omega$ such that $\omega x = \omega y$).
In 1977, Frankl and Deza \cite{ref:FranklDeza77} proved that if $S$ is an intersecting family of permutations of $\Omega$, then $|S| \leq (n-1)!$, with equality if and only if $S$ consists of all permutations sending a given point $\alpha$ to a given point $\beta$ \cite{ref:CameronKu03, ref:LaroseMalvenuto04}.

The first connection with derangements is direct.
If $S$ is a \emph{subgroup} of $\Sym(\Omega)$, then $S$ is an intersecting family of permutations if and only if $S$ has no derangements (exercise), so, for example, if $n > 1$, then such a subgroup is intransitive by Theorem~\ref{thm:existence}.

The second connection is indirect.
For $G \leq \Sym(\Omega)$, a point stabiliser $G_\omega$ of maximal size is an intersecting family of permutations inside $G$, but it need not be the largest one, as the following example demonstrates.

\begin{ex}
\label{ex:intersecting}
{\rm
Consider $G = A_4$ acting on the set $\Omega$ of six $2$-subsets of $\{1,2,3,4\}$, as in Examples~\ref{ex:k-subsets} and~\ref{ex:jehne}.
Then $|G_\omega| = 2$, but the subgroup $V_4$ of Example~\ref{ex:k-subsets} has no derangements, so it is an intersecting family of size $4$.
}
\end{ex}

To measure this failure, Li, Song and Pantangi \cite{ref:LiSongPantangi20} defined the \emph{intersection density} of a transitive group $G$ as
\[
\rho(G) = \max\{ |S|/|G_\omega| \mid \text{$S \subseteq G$ is an intersecting family} \}.
\]
How can we obtain upper bounds on $\rho(G)$?

Assume that $G \leq \Sym(\Omega)$ is transitive.
The \emph{derangement graph} of $G$, denoted $\Gamma_G$, is the graph whose vertex set is $G$ and where $x, y \in G$ are adjacent if $xy^{-1}$ is a derangement.
Then a subset $S$ of $G$ is an intersecting family if and only if it is a coclique of $\Gamma_G$.
Moreover, $\Gamma_G$ is the Cayley graph of $G$ with respect to the set of derangements, so $G$ acts vertex-transitively on it by right multiplication.
The Clique--Coclique Bound states that $|C| \cdot |S| \leq |G|$ for every clique $C$ and coclique $S$ of the vertex-transitive graph $\Gamma_G$.
Together with the Orbit--Stabiliser Theorem, this gives $\rho(G) \leq n/|C|$ for every clique $C$ of $\Gamma_G$.
Theorem~\ref{thm:existence} ensures that if $n > 1$, then $G$ has a derangement, so $\Gamma_G$ contains an edge and $\rho(G) \leq n/2$.
Moreover, in 2021, Meagher, Razafimahatratra and Spiga \cite{ref:MeagherRazafimahatratraSpiga21} proved that if $n > 2$, then $\Gamma_G$ contains a triangle, so $\rho(G) \leq n/3$, and they put forward the following conjecture as a question.

\begin{conj}
\label{conj:derangement_graph}
For all $k$, there exists $N_k$ such that if $G \leq \Sym(n)$ is transitive with $n > N_k$, then $\Gamma_G$ has a clique of size $k$.
\end{conj}

As just explained, Conjecture~\ref{conj:derangement_graph} implies the following.

\begin{conj}
\label{conj:intersecting}
For all $k$, there exists $N_k$ such that if $G \leq \Sym(n)$ is transitive with $n > N_k$, then $\rho(G) \leq n/k$.
\end{conj}

In 2026, Cazzola, Gogniat and Spiga \cite{ref:CazzolaGogniatSpiga26} proved these conjectures for $k = 4$ with $N_4 = 30$, which is best possible.

Let us now return to Conjecture~\ref{conj:neumann_praeger} arising from Neumann and Praeger's work on Kronecker equivalence.
Considering the formulation in (c), let $\overline{A}$ be the permutation group induced by $A$ on the cosets of $V$.
Every element of $U$ is $A$-conjugate into the point stabiliser $V$, so the image $\overline{U}$ of $U$, which has index $n$ in $\overline{A}$, has no derangements.
Hence the elements of a clique in $\Gamma_{\overline{A}}$ lie in distinct cosets of $\overline{U}$, so $|A:V| \leq N_{n+1}$ if Conjecture~\ref{conj:derangement_graph} holds.
Therefore, Conjecture~\ref{conj:derangement_graph} implies Conjecture~\ref{conj:neumann_praeger}, so Cazzola, Gogniat and Spiga \cite{ref:CazzolaGogniatSpiga26} deduced $h(3) = 10$ satisfies Conjecture~\ref{conj:neumann_praeger}(c) from $N_4 = 30$ satisfying Conjecture~\ref{conj:derangement_graph}.
Remarkably, in 2026, Anzanello and Spiga \cite{ref:AnzanelloSpiga} proved the converse, using the CFSG, so Conjectures~\ref{conj:derangement_graph} and~\ref{conj:neumann_praeger} are equivalent.

\subsection{Game theory: Isbell's Conjecture}
\label{ss:connections_isbell}

\sectionquestion{You ask a room of people whether to buy red or white wine. How can you fairly decide what to buy?}

Fix a positive integer $n$ and let $\Omega = \{ 1, \dots, n \}$.
A \emph{set of winning coalitions} on $\Omega$ is a subset $W \subseteq \P(\Omega)$ such that for all $A, B \in \P(\Omega)$
\begin{itemize}\setlength{\itemsep}{-4pt}
\item[{\rm (a)}]
$A \in W$ if and only if $\Omega \setminus A \not\in W$
\item[{\rm (b)}]
if $A \in W$ and $A \subseteq B$, then $B \in W$. 
\end{itemize}
The \emph{degree} of $W$ is $n$. 
(A set of winning coalitions is exactly a \emph{maximal} intersecting family of subsets of $\Omega$, see Section~\ref{ss:connections_intersecting}.)

Returning to our problem, we number the people $1$ to $n$ and assume that we have fixed a set $W$ of winning coalitions of degree~$n$.
The first condition ensures that the set of people who voted ``red'' will be in $W$ or (its complement) the set of people who voted ``white'' will be in $W$, but not both; we choose the option that is contained in $W$. 
The second condition means that if red wins and people switch from white to red, then red still wins.

\begin{ex}
\label{ex:isbell_dictatorship}
{\rm
For any $n$, an example of a set of winning coalitions of degree~$n$ is a \emph{dictatorship}, namely for fixed $i \in \Omega$ the set 
\[
W = \{ A \in \P(\Omega) \mid i \in A\}.
\]
Here, whatever the dictator $i$ voted for is chosen.
}
\end{ex}

A set of winning coalitions is \emph{fair} if the group $G \leq \Sym(n)$ that stabilises $W$, namely
\[
\Aut(W) = \{ g \in \Sym(n) \mid \text{$Ag \in W$ for all $A \in W$}\}
\]
acts transitively on $\Omega$.
A dictatorship is not fair: the dictator is a fixed point of the group stabilising $W$.

\begin{ex}
\label{ex:isbell_majority}
{\rm
For odd $n$, an example of a fair set of winning coalitions of degree~$n$ is \emph{majority rules}, namely 
\[
W = \{ A \in \P(\Omega) \mid |A| > n/2 \}.
\]
Here, whichever option gets the most votes is chosen.
}
\end{ex}

When $n$ is even, there need not exist a fair set of winning coalitions at all.
Indeed, there is none when $n = 4$ (exercise).

Sets of winning coalitions are the \emph{simple games} of von~Neumann and Morgenstern \cite{ref:vonNeumannMorgenstern44}.
Isbell studied these in the 1950s and 1960s \cite{ref:Isbell57, ref:Isbell60, ref:Isbell64}, and he made the following conjecture \cite{ref:Isbell60}.

\begin{conj}
\label{conj:isbell_coalitions}
There is a map $m \colon \Nat\to \Nat$ such that there is no fair set of winning coalitions of degree $2^a b$ with $b$ odd and $a \geq m(b)$.
\end{conj}

The resemblance between Conjectures~\ref{conj:isbell} and~\ref{conj:isbell_coalitions} of Isbell is no coincidence (see \cite[Lemma~1]{ref:Isbell57}).

\begin{lem}
\label{lem:isbell}
Let $n$ be a positive integer.
Then there exists a fair set of winning coalitions of degree $n$ if and only if there exists a transitive group $G \leq \Sym(n)$ with no derangement of $2$-power order.
\end{lem}

\begin{proof}
First assume that there exists a fair set of winning coalitions $W$ of degree $n$.
Let $G = \Aut(W)$.
Since $W$ is fair, $G \leq \Sym(n)$ is transitive.
Suppose that $G$ has a derangement $g$ of $2$-power order.
Then every cycle of $g$ has even length, so (by making some choices) we can fix a set $A \in \P(\Omega)$ that contains every second element from every cycle of $g$. 
Then $Ag = \Omega \setminus A$, and since $g \in \Aut(W)$, this forces $\Omega \setminus A \in W$.
Therefore, either $A, \Omega \setminus A \in W$ or $A, \Omega \setminus A \not\in W$, but $W$ is a set of winning coalitions, so exactly one of $A$ and $\Omega \setminus A$ is contained in $W$, which is a contradiction.
Therefore, $G$ does not contain a derangement of $2$-power order.

Conversely, suppose that $G \leq \Sym(n)$ is transitive with no derangement of $2$-power order.
Let 
\[
W_0 = \{ A \in \P(\Omega) \mid |A| > n/2 \}.
\]
If $n$ is odd, then $W = W_0$ is a fair set of winning coalitions.
We now assume that $n$ is even.
Write $X = \{ A \in \P(\Omega) \mid |A| = n/2 \}$.

We claim that $Ag \neq \Omega \setminus A$ for all $A \in X$ and $g \in G$.
To see this let $A \in X$ and $g \in G$ and suppose that $Ag = \Omega \setminus A$.
Then $g$ is a derangement and every cycle of $g$ has even length, so raising $g$ to the odd part of its order gives a derangement of $2$-power order (see Remark~\ref{rem:p-elements}), which is a contradiction.

With this claim in place, by making some choices, we can write the set of orbits of $G$ on $X$ as the disjoint union 
\[
\{ \Delta_1, \dots, \Delta_k \} \cup \{ \Gamma_1, \dots, \Gamma_k \},
\]
where $\Gamma_i = \{ \Omega \setminus A \mid A \in \Delta_i \}$ for each $i$.
Define $W = W_0 \cup \bigcup_{i=1}^k \Delta_i$.

We claim that $W$ is a fair set of winning coalitions.
For the first condition, if $A \in W$, then $A \in \Delta_i$ for some $i$ or $|A| > n/2$, so $\Omega \setminus A \in \Gamma_i$ or $|A| < n/2$; in either case, $\Omega \setminus A \not\in W$.
Conversely, if $\Omega \setminus A \not\in W$, then $\Omega \setminus A \in \Gamma_i$ for some $i$ or $|\Omega \setminus A| < n/2$, so $A \in \Delta_i$ or $|A| > n/2$; in either case, $A \in W$.
For the second condition, if $A \in W$ and $A \subseteq B$, then $A = B$ or $|B| > |A| \geq n/2$; in either case, $B \in W$.
For the fairness condition, note that $G$ stabilises $W$, so $G \leq \Aut(W)$; since $G$ is transitive, so is $\Aut(W)$.
\end{proof}

We conclude with an example to highlight that it is possible to have a fair set of winning coalitions of even degree.

\begin{ex}
\label{ex:isbell_6}
{\rm
By Example~\ref{ex:isbell}, if $q \equiv 1 \pmod{4}$, then $G = \PSL_2(q)$ has no derangements of $2$-power order on the set $\Omega$ of $1$-spaces of $\F_q^2$.
Therefore, Lemma~\ref{lem:isbell} gives a fair set of winning coalitions of degree $|\Omega| = q + 1$.

We give the details for $q = 5$.
View $\Omega$ as $\F_5 \cup \{\infty\}$, as in Example~\ref{ex:pgl2}.
For distinct $a, b, c \in \F_5$, the set $\{ a, b, c \}$ is square if $(a - b)(b - c)(c - a)$ is a square in $\F_5^\times$, and $\{ a, b, \infty \}$ is square if $a - b$ is a square in $\F_5^\times$.
Since $-1$ is a square in $\F_5^\times$, whether a $3$-subset is square is independent of the order the elements are given in.
We claim that 
\[
W = \{ A \in \P(\Omega) \mid \text{$|A| > 3$ or $|A| = 3$  and $A$ is square}\}
\] 
is a fair set of winning coalitions on $\Omega$.
Since $G$ has no derangements of $2$-power order, no element of $G$ maps a set to its complement, so $G$ has at least two orbits on $3$-subsets of $\Omega$ and, in fact, exactly two since $G$ has index two in $\PGL_2(5)$, which is (sharply) $3$-transitive on $\Omega$.
Moreover, being square distinguishes these two orbits (exercise).
Accordingly, $W$ is the set of winning coalitions constructed in the proof of Lemma~\ref{lem:isbell}.
}
\end{ex}

\subsection{Graph theory: the Polycirculant Conjecture}
\label{ss:connections_polycirculant}

\sectionquestion{Can every finite vertex-transitive graph be drawn with a rotational symmetry that moves every vertex?}

A \emph{polycirculant} is a graph with a nontrivial automorphism all of whose cycles have the same length.
Drawing these cycles as concentric regular polygons turns this automorphism into a rotation, as for the Petersen graph and the automorphism $g = (1 \ 2 \ 3 \ 4 \ 5)(6 \ 7 \ 8 \ 9 \ 10)$ in Figure~\ref{fig:petersen}.
In contrast, a path with an odd number of vertices (at least three) is not a polycirculant, since its only nontrivial automorphism fixes the central vertex.

One appealing feature of polycirculants is that they admit an efficient description.
Let $\Gamma$ be a polycirculant with an automorphism consisting of $r$ cycles of length $m$.
The \emph{Frucht notation} \cite{ref:Frucht70} for $\Gamma$ is the following labelled graph.
For each cycle $(v_1 \ \dots \ v_m)$, there is a vertex labelled $m \mid k$ where $k$ is chosen such that $v_i$ is adjacent to $v_{i+k}$ for all $i$ (where indices are always interpreted modulo $m$).
For distinct cycles $(u_1 \ \dots \ u_m)$ and $(v_1 \ \dots \ v_m)$ there is an edge labelled $k$ if there is an edge between $u_i$ and $v_{i+k}$ for all $i$, and no edge otherwise (in which case it follows that none of $u_1, \dots, u_m$ are adjacent to any of $v_1, \dots, v_m$); labels of $0$ are omitted from edges.
The Frucht notation for the Petersen graph for the automorphism $g$ is given in Figure~\ref{fig:petersen}.

In 1981, Maru\v{s}i\v{c}~\cite{ref:Marusic81} asked the following (actually for all digraphs), which was later asked again by Leighton~\cite{ref:Leighton83} and Jordan~\cite{ref:Jordan88}.
Here, a graph $\Gamma$ is \emph{vertex-transitive} if the automorphism group $\Aut(\Gamma)$ acts transitively on the set of vertices of $\Gamma$.

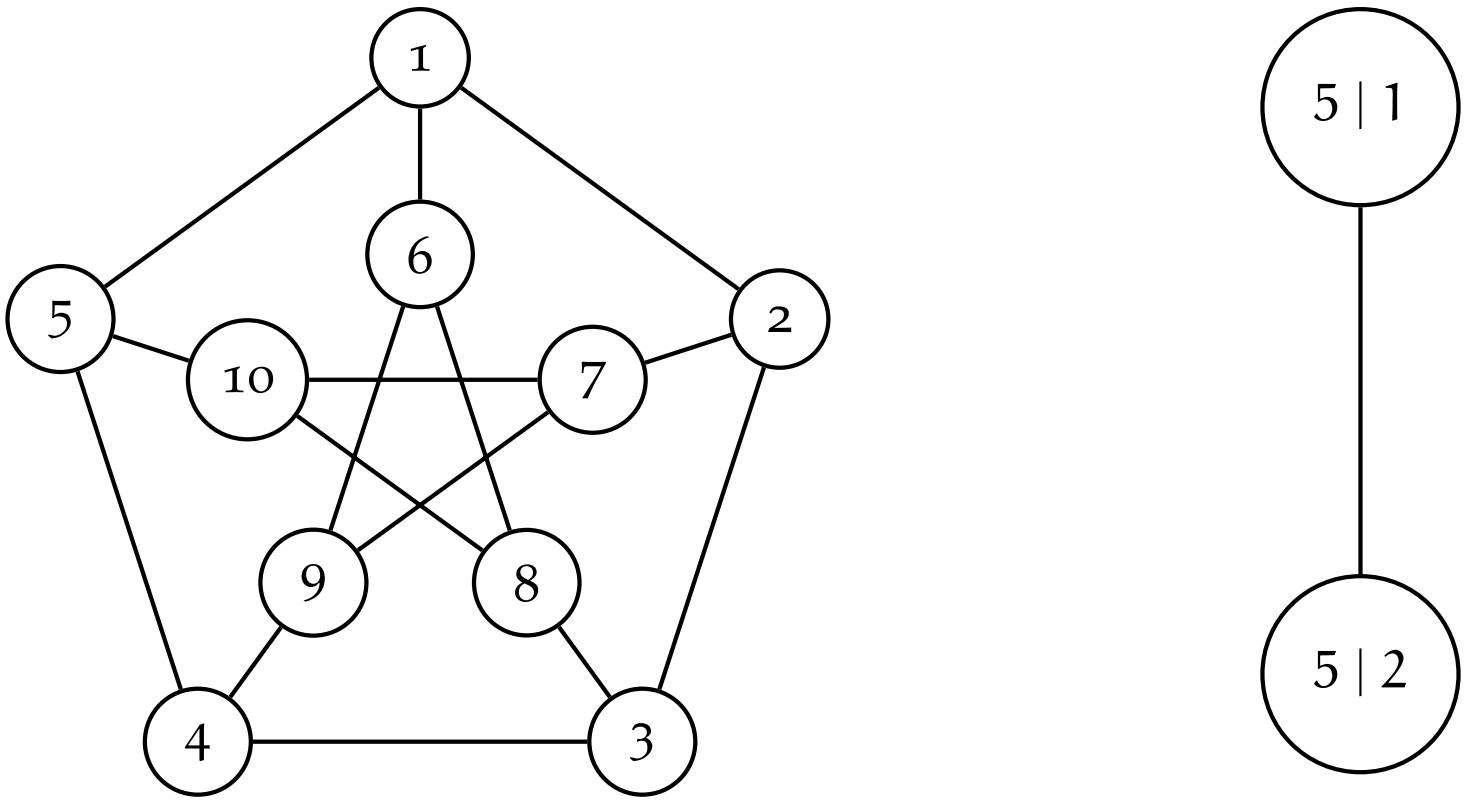
\begin{figure}[t]
\begin{multicols}{2}
\centering
\begin{tikzpicture}[every node/.style={circle, draw, fill=white, inner sep=4pt, font=\small}, thick]
\def\R{2.5}
\def\r{1.2}
\node (1) at ({90}:\R)     {1};
\node (2) at ({90-72}:\R)  {2};
\node (3) at ({90-144}:\R) {3};
\node (4) at ({90-216}:\R) {4};
\node (5) at ({90-288}:\R) {5};
\node (6)  at ({90}:\r)     {6};
\node (7)  at ({90-72}:\r)  {7};
\node (8)  at ({90-144}:\r) {8};
\node (9)  at ({90-216}:\r) {9};
\node (10) at ({90-288}:\r) {10};
\draw (1)--(2)--(3)--(4)--(5)--(1);
\draw (6)--(8)--(10)--(7)--(9)--(6);
\draw (1)--(6);
\draw (2)--(7);
\draw (3)--(8);
\draw (4)--(9);
\draw (5)--(10);
\end{tikzpicture}
\columnbreak

\begin{tikzpicture}[every node/.style={circle, draw, inner sep=6pt, font=\small},thick]
\node (top) at (0,3.75) {$5 \mid 1$};
\node (bot) at (0,0) {$5 \mid 2$};
\draw (top)--(bot);
\end{tikzpicture}
\end{multicols}
\caption{The Petersen graph and its Frucht notation}
\label{fig:petersen}
\end{figure}

\begin{que}
\label{que:polycirculant}
Is every finite vertex-transitive graph a polycirculant?
\end{que}

Question~\ref{que:polycirculant} has been answered affirmatively in many special cases, such as for graphs of valency at most~$4$ \cite{ref:DobsonMalnicMarusicNowitz07II, ref:GiudiciMorganPotocnikVerret15,ref:MarusicScapellato98}, for graphs with a square-free number of vertices \cite{ref:DobsonMalnicMarusicNowitz07I} and for distance-transitive graphs \cite{ref:KutnarSparl10}.

Let us interpret this problem group theoretically.

\begin{lem}
\label{lem:polycirculant}
A finite graph is a polycirculant if and only if its automorphism group has a derangement of prime order.
\end{lem}

\begin{proof}
Let $\Gamma$ be a finite graph.
If $\Aut(\Gamma)$ has a derangement $g$ of prime order $p$, then every cycle of $g$ has length $p > 1$, so $\Gamma$ is a polycirculant.
Conversely, if $\Gamma$ is a polycirculant and every cycle of $g \in \Aut(\Gamma)$ has length $r > 1$, then let $p$ be a prime divisor of $r$ and note that every cycle of $g^{r/p}$ has length $p$, so $g^{r/p}$ is a derangement of prime order.
\end{proof}

In Section~\ref{ss:derangements_prime}, we saw that there are transitive permutation groups with no derangements of prime order.
However, this does not immediately give a negative answer to Question~\ref{que:polycirculant}, since not all transitive permutation groups arise as automorphism groups of graphs, as the following example highlights.

\begin{ex}
\label{ex:polycirculant}
{\rm
The primitive group $G = {\rm M}_{11} \leq \Sym(12)$ has no derangements of prime order (see (3d) in the proof of Theorem~\ref{thm:fein_kantor_schacher}).
Suppose that $\Gamma$ were a graph on $12$ vertices with $\Aut(\Gamma) = G$. 
Since $G$ is $2$-transitive (in fact, it is $3$-transitive), given any two pairs of distinct vertices $(u_1, u_2)$ and $(v_1, v_2)$ of $\Gamma$, there is an automorphism $g \in G$ such that $(u_1, u_2)g = (v_1, v_2)$.
However, $g$ is an automorphism, so either there is an edge between all pairs of distinct vertices of $\Gamma$ (so $\Gamma$ is complete) or between no pairs of distinct vertices of $\Gamma$ (so $\Gamma$ is empty).
In either case, $\Aut(\Gamma) = \Sym(12)$, which is a contradiction.
}
\end{ex}

To generalise Example~\ref{ex:polycirculant}, define the \emph{$2$-closure} of $G \leq \Sym(n)$, denoted $G^{(2)}$, to be the largest subgroup of $\Sym(n)$ that stabilises the orbitals of $G$, and say that $G$ is \emph{$2$-closed} if $G^{(2)} = G$.
For example, if $G \leq \Sym(n)$ is $2$-transitive, then $G^{(2)} = \Sym(n)$.
The automorphism group of a graph is $2$-closed as its edge set is a union of orbitals.
The converse fails, but every $2$-closed group is the automorphism group of an edge-coloured digraph.

In 1997, Klin conjectured the following (see \cite{ref:Cameron97}), generalising Maru\v{s}i\v{c}'s question above.

\begin{conj}[Polycirculant Conjecture]
\label{conj:polycirculant}
Every transitive $2$-closed group $G \leq \Sym(n)$ with $n > 1$ has a derangement of prime order.
\end{conj}

Progress is surveyed in \cite[Section~4]{ref:ArezoomandAbdollahiSpiga19} and \cite[Section~1.3.5]{ref:BurnessGiudici16}, but in July 2026, shortly after the minicourse, Freedman and Lee \cite{ref:FreedmanLee} gave a negative answer to Maru\v{s}i\v{c}'s question and hence disproved the Polycirculant Conjecture.
Building on \cite{ref:ChenLeeMitrovicOBrienXia26}, they give an infinite family of examples, the smallest of which is the following.

\begin{ex}
\label{ex:polycirculant_counterexample}
{\rm
This is a variant of Example~\ref{ex:elusive_odd}.
We build on the degree $336$ action of $\PSU_3(3)$ on $\PSU_3(3)/K$, where $K \cong C_3 \times S_3$, in which every derangement of prime order has order~$7$.
Let $N = \F_7^6$ be the unique irreducible $6$-dimensional module for $\PSU_3(3)$ over $\F_7$.
In \cite{ref:FreedmanLee}, a nonsplit extension $G = N.\PSU_3(3)$ is considered in which every element of order $7$ lies in $N$, and a $4$-dimensional $K$-submodule $N_0$ of $N$ is given whose $\PSU_3(3)$-translates cover $N$ (both facts are verified computationally).
Let $H = N_0{:}K$ and consider the degree $16464 = 2^4 \cdot 3 \cdot 7^3$ action of $G$ on $\Omega = G/H$.
An element of order $7$ lies in $N$, hence in a conjugate of $N_0 \leq H$, while an element of order $2$ or $3$ is conjugate into $N{:}K$, since this holds for its image in $\PSU_3(3)$, and hence into $K \leq H$ by the Schur--Zassenhaus Theorem.
Hence, $G$ is elusive.
A computation with the orbitals of $G$ shows that $G$ is the automorphism group of a graph $\Gamma$ of valency~$63$, so $G$ is $2$-closed, $\Gamma$ gives a negative answer to Maru\v{s}i\v{c}'s question and $G$ is a counterexample to the Polycirculant Conjecture.
Unlike in Example~\ref{ex:elusive_odd}, $N_0$ has codimension $2$ rather than $1$, which is necessary for a $2$-closed elusive group by \cite[Theorem~5.3]{ref:CameronGiudiciJonesKantorKlinMarusicNowitz02} (see \cite[Remark~1.5]{ref:FreedmanLee}).
}
\end{ex}

\clearpage

\section*{Acknowledgements}

The author thanks the organisers of GABY 2026 for the invitation to write this survey, the opportunity to deliver the minicourse that it is based on and their hospitality during the summer school itself.
He also thanks Jay Taylor for conversations surrounding Section~\ref{ss:connections_character_restriction} and Marco Barbieri and Julian Feuerpfeil for their helpful comments on an earlier draft.
The author used Anthropic's Claude (Opus~5.5 and Fable~5.1) to help him understand some of the ideas in this survey, to help him find some of the works referenced in this survey, to help him write this survey and to review a draft of the survey.


\renewcommand{\bibsection}{

\vspace{11pt}

\noindent 
Scott Harper \\
School of Mathematics \\
University of Birmingham \\
Birmingham, B15 2TT (United Kingdom) \\
s.harper.3@bham.ac.uk 

}
\end{document}